\documentclass[a4paper,12pt,leqno]{amsart}
\usepackage{amsmath,amssymb,mathrsfs}
\usepackage[matrix,arrow]{xy} 
\usepackage{nameref}
\usepackage{xcolor,hyperref}
\hypersetup{colorlinks=true,citecolor=black,linkcolor=black}
\newcommand{\h}{\hbox}
\newcommand{\q}{\quad}
\newcommand{\nin}{\par\noindent}

\newcommand{\ms}{\par\medskip}
\newcommand{\sk}{\par\smallskip}

\newcommand{\msn}{\par\medskip\noindent}

\newcommand{\mcup}{\h{$\bigcup$}}
\newcommand{\msum}{\h{$\sum$}}
\newcommand{\mprod}{\h{$\prod$}}
\newcommand{\mopl}{\h{$\bigoplus$}}
\newcommand{\ssb}{\raise.15ex\h{${\scriptscriptstyle\bullet}$}}
\newcommand{\ssc}{\,\raise.15ex\h{${\scriptstyle\circ}$}\,}
\newcommand{\ds}{\rlap{\raise.4ex\h{$\downarrow$}}{\downarrow}}
\newcommand{\ee}{{\mathbf e}}
\newcommand{\C}{{\mathbb C}}
\newcommand{\N}{{\mathbb N}}
\newcommand{\PP}{{\mathbb P}}
\newcommand{\Q}{{\mathbb Q}}

\newcommand{\Rb}{{\mathbf R}}
\newcommand{\Z}{{\mathbb Z}}
\newcommand{\B}{{\mathcal B}}
\newcommand{\D}{{\mathcal D}}
\newcommand{\F}{{\mathcal F}}
\newcommand{\G}{{\mathcal G}}
\newcommand{\Hc}{{\mathcal H}}
\newcommand{\M}{{\mathcal M}}
\newcommand{\Lc}{{\mathcal L}}
\newcommand{\OO}{{\mathcal O}}
\newcommand{\R}{{\mathcal R}}
\newcommand{\Sc}{{\mathcal S}}

\newcommand{\alt}{\widetilde{\alpha}}

\newcommand{\Ct}{\widetilde{C}}

\newcommand{\Et}{\widetilde{E}}
\newcommand{\Gt}{\widetilde{\mathcal G}}
\newcommand{\Ht}{\widetilde{H}}
\newcommand{\Kt}{\widetilde{K}}

\newcommand{\Mt}{\widetilde{M}}
\newcommand{\Nt}{\widetilde{N}}

\newcommand{\Xt}{\widetilde{X}}

\newcommand{\f}{\tilde{f}}

\newcommand{\dd}{{\partial}}
\newcommand{\ddd}{{\rm d}}

\newcommand{\dfw}{{\rm d}f{\wedge}}
\newcommand{\dti}{\partial_t^{-1}}
\newcommand{\spKf}{{}^{s}{}'\!K_{f'}^{\ssb}}
\newcommand{\kod}{\tfrac{k}{d}}
\newcommand{\nod}{\tfrac{n}{d}}
\newcommand{\ood}{\tfrac{1}{d}}
\newcommand{\al}{\alpha}
\newcommand{\de}{\delta}
\newcommand{\ga}{\gamma}

\newcommand{\la}{\lambda}
\newcommand{\om}{\omega}
\newcommand{\Om}{\Omega}
\newcommand{\Si}{\Sigma}
\newcommand{\Gr}{{\rm Gr}}
\newcommand{\Sp}{{\rm Sp}}
\newcommand{\Supp}{{\rm Supp}}
\newcommand{\Ff}{F_{\!f}}
\newcommand{\Fhz}{F_{h_z}}
\newcommand{\stm}{\,{\setminus}\,}
\newcommand{\sst}{\,{\subset}\,}
\newcommand{\defs}{\,{:=}\,}
\newcommand{\nes}{\,{\ne}\,}
\newcommand{\eq}{\,{=}\,}
\newcommand{\ins}{\,{\in}\,}
\newcommand{\pl}{\one {+}\one}
\newcommand{\mi}{\one {-}\one}
\newcommand{\bl}{\bigl}
\newcommand{\br}{\bigr}
\newcommand{\into}{\hookrightarrow}
\newcommand{\onto}{\twoheadrightarrow}
\newcommand{\simto}{\,\,\rlap{\hskip1.5mm\raise1.4mm\hbox{$\sim$}}\hbox{$\longrightarrow$}\,\,}
\newcommand{\ges}{\geqslant}
\newcommand{\les}{\leqslant}
\newcommand{\gess}{\,{\geqslant}\,}
\newcommand{\less}{\,{\leqslant}\,}
\newcommand{\sgt}{\,{>}\,}
\newcommand{\slt}{\,{<}\,}
\newcommand{\one}{\hskip1pt}
\newcommand{\mm}{{\mathfrak m}}

\makeatletter
\renewcommand\section{\@startsection{section}{1}{0pt}{-3.5ex plus -1ex minus -.2ex}{-2.3ex plus -.2ex}{\normalfont\bfseries}}
\makeatother

\theoremstyle{plain}
\newtheorem{thm}{Theorem}[section]

\newtheorem{prop}{Proposition}[section]
\newtheorem{lem}{Lemma}[section]
\newtheorem{ithm}{Theorem}
\newtheorem{icor}{Corollary}
\newtheorem{iprop}{Proposition}

\theoremstyle{definition}
\newtheorem{rem}{Remark}[section]

\newtheorem{irem}{Remark}

\newtheorem{A}{A}

\begin{document}
\title[Bernstein-Sato polynomials]{Bernstein-Sato polynomials\\for projective hypersurfaces with\\weighted homogeneous isolated singularities}
\author{Morihiko Saito}
\address{RIMS Kyoto University, Kyoto 606-8502 Japan}
\begin{abstract} We present a quite efficient method to calculate the roots of Bernstein-Sato polynomial for a defining polynomial $f$ of a projective hypersurface $Z\subset{\mathbb P}^{n-1}$ of degree $d$ having only weighted homogeneous isolated singularities. We prove the $E_2$-degeneration of the pole order spectral sequence so that the computation of roots is reduced to the one of the Hilbert series of the Jacobian ring of $f$ except the special case where $f$ is annihilated by a nonzero vector field on ${\mathbb C}^n$ with linear function coefficients. In the three variable case with $d>4$ we may assume that this vector field is a linear combination of $x\partial_x, y\partial_y, z\partial_z$, where $f$ is called extremely degenerated; in particular, the latter case does not contain any essential indecomposable central hyperplane arrangement in ${\mathbb C}^3$. Combined with the self-duality of the Koszul complex and a theorem of Dimca and Popescu, it implies for $n=3$ with $d>4$ except the extremely degenerated case that $R_f=\frac{1}{d}({\mathbb Z}\cap[3,k'])\cup R_Z$. Here $R_f,R_Z$ are the roots of Bernstein-Sato polynomials of $f$ and $Z$ up to sign, and $k'\eq\max(2d-3,k_{\max}+3)$ with $k_{\max}$ the maximal degree of the ``torsion part" of the Jacobian ring, where the latter is known to be at most $2d-5$ in the hyperplane arrangement case.
\end{abstract}
\maketitle

\part*{Introduction} \label{intr} \nin
Let $Z$ be a projective hypersurface in $\PP^{n-1}$ defined by a homogeneous polynomial $f$ of $n$ variables with $n\gess 3$ and $d\defs\deg f\gess 3$. We assume that $Z$ is reduced and $f$ is essential (that is, $f$ is not a polynomial of $n{-}1$ variables). Let $b_f(s)$ be the Bernstein-Sato polynomial of $f$. Set
$$\R_f\defs\bl\{\al\in\Q\mid b_f(-\al)=0\br\}\subset\Q_{>0},$$
and similarly for $\R_{h_z}$ replacing $f$ by a local defining function $h_z$ of $(Z,z)\subset(Y,z)$ ($z\in Z$). Define the set of roots of Bernstein-Sato polynomial {\it supported at the origin} (up to sign) by
\begin{equation} \label{1}
\R_f^0:=\R_f\setminus\R_Z\,\q\h{with}\q\R_Z:=\mcup_{z\in{\rm Sing}\,Z}\,\R_{h_z}\,\subset\,\R_f.
\end{equation}
Here the last inclusion follows from the equality $b_{h_z}(s)=b_{f,y}(s)$ with $b_{f,y}(s)$ the local Bernstein-Sato polynomial of $f$ at $y\in\C^n\setminus\{0\}$ with $[y]=z$ in $Y=(\C^n\setminus\{0\})/\C^*$. (This follows from the assertion that $b_{f,y}(s)$ depends only on $f^{-1}(0)$, see Remark~\ref{R20.1} below.)
\sk
Set $\Ff:=f^{-1}(1)\subset X\defs\C^n$, the Milnor fiber of $f$. We have the pole order filtration $P$ on each monodromy eigenspace $H^j(\Ff,\C)_{\la}\defs{\rm Ker}(T_s\mi\la)\sst H^j(\Ff,\C)$ ($\la\in\C^*$) with $T\eq T_sT_u$ the Jordan decomposition of the monodromy. Recall the following.

\begin{ithm}[{\cite[Theorem 2]{bcm}}] \label{T1}
Assume $\al\notin\R_Z$. If $\al$ satisfies the condition
\begin{equation} \label{2}
\al\notin\R_Z+\Z_{<0},
\end{equation}
then
\begin{equation} \label{3}
\al\in\R_f^0\iff\Gr_P^p\,H^{n-1}(\Ff,\C)_{\la}\ne 0\,\,\,\bl(p=[n{-}\al],\,\,\la=e^{-2\pi i\al}\br).
\end{equation}
If condition~{\rm\eqref{2}} does not hold, then only the implication $\Longleftarrow$ holds in {\rm\eqref{3}}.
\end{ithm}

\begin{irem} \label{R1}
There are examples of {\it non-reduced\one} hyperplane arrangements of 3 variables such that condition~\eqref{2} is unsatisfied and the implication $\Longrightarrow$ in \eqref{3} fails, see \cite[Example 4.5]{nwh}. This problem is related to the {\it asymptotic expansion} of $f^s\in\D_X[s]f^s$ via the $V$-filtration, and is rather complicated, see Remark~\ref{R20.3} below.
\end{irem}

Theorem~\ref{T1} gives a partial generalization of a well-known theorem in the isolated singularity case asserting that the Steenbrink spectral numbers coincide with the roots of the microlocal Bernstein-Sato polynomial $\widetilde{b}_f(s):=b_f(s)/(s{+}1)$ up to a sign (forgetting the multiplicities) in the case $f$ is a weighted homogeneous polynomial, see \eqref{7.2} below.
\sk
We have the {\it pole order spectral sequence} associated with the {\it pole order filtration} on the (algebraic) {\it microlocal Gauss-Manin complex} of $f$. Set
$$R:=\C[x_1,\dots,x_n],$$
with $x_1,\dots,x_n$ the coordinates of $\C^n$.
This is a graded ring with $\deg x_i\eq1$.
Let $\Om^{\ssb}$ be the complex of the exterior products of the K\"ahler differentials of $R$ over $\C$ so that $\Om^p$ is a free $R$-module of rank $\binom{n}{p}$.
This complex has anti-commuting two differentials $\ddd$ and $\dfw$ preserving the grading (up to the shift by $d=\deg f$ in the case of $\dfw$), where each component of the complex is graded with $\deg x_i=\deg\ddd x_i=1$.
\sk
In this paper we assume the {\it isolated singularity condition\,}:
\msn
(IS)\q\q\q\q\q\q\q\q\q$\sigma_Z:=\dim\Si=0\q\h{with}\q\Si:={\rm Sing}\,Z$.
\msn
We have the vanishing
\begin{equation} \label{4}
H^j_{\ddd}\bl(H^{\ssb}_{\dfw}\Om^{\ssb}\br)=0\q\h{unless}\,\,\,\,j=n{-}1\,\,\,\,\h{or}\,\,\,\,n,
\end{equation}
where $H^{\ssb}_{\dfw}$ means that the cohomology with respect to the differential $\dfw$ is taken, and similarly for $H^{\ssb}_{\ddd}$. These are members of the $E_2$-term of the spectral sequence associated with the double complex with anti-commuting two differentials $\ddd$ and $\dfw$ on $\Om^{\ssb}$, see also \eqref{1.1} below. The latter has been studied in \cite{Di1}, although the relation to the Gauss-Manin system and the Brieskorn modules was not mentioned there. Note that the {\it usual} (that is, non-microlocal) pole order spectral sequence was considered there, but its $E_2$-degeneration is equivalent to that of the {\it microlocal} one (see \cite[Corollary 4.7]{kosz}), and the latter is related to the pole order filtration as in \eqref{2.8} below. It has been observed in many examples by A.~Dimca and G.~Sticlaru \cite{DiSt2} that the pole order spectral sequence degenerates at $E_2$ if the singularities of $Z$ are isolated and weighted homogeneous. We have the following.

\begin{iprop} \label{P1}
Under the assumption {\rm (IS)}, the $E_2$-degeneration of the pole order spectral sequence is equivalent to each of the following two conditions\,$:$
\begin{equation} \label{5}
\dim H^n_{\ddd}\bl(H^{\ssb}_{\dfw}\Om^{\ssb}\br)=\dim H^{n-1}(\Ff,\C).
\end{equation}
\vskip-6mm
\begin{equation} \label{6}
\dim H^{n-1}_{\ddd}\bl(H^{\ssb}_{\dfw}\Om^{\ssb}\br)=\dim H^{n-2}(\Ff,\C).
\end{equation}
If these equivalent conditions are satisfied, then the pole order spectrum $\Sp_P(f)$ {$($see \cite{kosz}$)$} is given by the difference of the Hilbert series of the graded $\C$-vector spaces
$$M_{\ssb}^{(2)}:=H^n_{\ddd}\bl(H^{\ssb}_{\dfw}\Om^{\ssb}\br),\q N_{\ssb}^{(2)}:=H^{n-1}_{\ddd}\bl(H^{\ssb}_{\dfw}\Om^{\ssb}\br)(-d),$$
with variable $t$ of the Hilbert series replaced by $t^{1/d}$ {$($}where $(m)$ for $m\in\Z$ denotes the shift of grading so that $E(m)_k:=E_{m+k}$ for any graded module $E)$, and moreover, in the notation of {\rm\eqref{3}}, we have the canonical isomorphisms
\begin{equation} \label{7}
M^{(2)}_k=\Gr_P^p\,H^{n-1}(\Ff,\C)_{\la}\q\h{for}\,\,\,\al=\kod,\,\,p=[n-\al],\,\,\la=e^{-2\pi i\al}.
\end{equation}
\end{iprop}

In this paper we prove the following.

\begin{ithm} \label{T2}
Under the assumption {\rm (IS)}, assume further
\msn
\rlap{\rm(WH)}\hskip1.5cm\hangindent=1.5cm\hangafter=1
Every singularity of $Z$ is analytically defined by a weighted homogeneous polynomial\\
with an isolated singularity, see Remark~{\rm\ref{R4.1}} below.
\msn
Then the pole order spectral sequence degenerates at $E_2$ so that $\eqref{7}$ holds.
\end{ithm}

This is quite nontrivial even in the case $n\eq3$. Without assuming (WH), Theorem~\ref{T2} never holds, see \cite[Theorem~5.2]{kosz}. (Here ``weighted homogeneity" cannot be replaced by ``semi-weighted-homogeneity".) Recall that the assumption~(WH) is equivalent to that every singularity of $Z$ is {\it quasihomogeneous}, that is, $h_z\in(\dd h_z)$ ($\forall\,z\in\Si$), where $h_z$ is as in \eqref{1}, and $(\dd h_z)$ is the Jacobian ideal generated by the partial derivatives of $h_z$, see \cite{SaK} (and Remark~\ref{R4.1} below). Note that $\R_{h_z}$ can be determined only by the {\it weights} in the weighted homogeneous isolated singularity case, see \eqref{7.2}--\eqref{7.3} and A.\ref{A.1} in Appendix below. This is {\it quite different\,} from the non-quasihomogeneous isolated hypersurface singularity case where we need a computer program to determine the local Bernstein-Sato polynomials of $Z$ in general.
\sk
By \cite[Corollary 4.7]{kosz}, Theorem~\ref{T2} implies the following.

\begin{icor} \label{C1}
Under the hypotheses~{\rm(IS)}-{\rm(WH)}, the Brieskorn module $H^nA_f^{\ssb}$ in the notation of \cite{kosz} is torsion-free.
\end{icor}

Combining Theorem~\ref{T2} with \cite{kosz}, \cite{bcm}, we can show the following.

\begin{icor} \label{C2}
Under the assumptions~{\rm(IS)}-{\rm(WH)}, any $\al\in\R^0_f$ satisfying condition~{\rm\eqref{2}} in Theorem~$\ref{T1}$ can be detected by using the Hilbert series of the graded $\C$-vector space $M^{(2)}$ in Proposition~$\ref{P1}$.
\end{icor}

Set
$$M_{\ssb}:=H^n_{\dfw}\Om^{\ssb},\q N_{\ssb}:=H^{n-1}_{\dfw}\Om^{\ssb}(-d).$$
These are graded $\C$-vector spaces. Let $y$ be a sufficiently general linear combination of the coordinates $x_1,\dots,x_n$ of $\C^n$, and $M'_{\ssb}\subset M_{\ssb}$ be the $y$-torsion part. Set $M''_{\ssb}:=M_{\ssb}/M'_{\ssb}$. Note that $M'_{\ssb}$ is a finite dimensional graded $\C$-vector subspace of $M_{\ssb}$, and is equal to the zeroth local cohomology $H^0_{\mm}M_{\ssb}$ with $\mm\sst R$ the maximal ideal at 0, see \cite{kosz}.

\begin{ithm} \label{T3}
Under the assumptions~{\rm(IS)}-{\rm(WH)}, we have the injectivity of the composition of canonical morphisms
$$M'_{\ssb}\into M_{\ssb}\onto M_{\ssb}^{(2)}.$$
\end{ithm}

Using Theorems~\ref{T1}, \ref{T2}, \ref{T3} together with \cite[Corollary 2 and Theorem 5.3]{kosz}, we can deduce

\begin{icor} \label{C3}
Under the hypotheses~{\rm(IS)}-{\rm(WH)}, assume $M'_k\ne 0$ for some $k\ins\Z$. Then $M^{(2)}_k\ne 0$, hence $\kod\in\R_f$. The converse holds if we assume $k>(n{-}1)d\mi n$ and moreover $\kod\notin\R_Z$ in the case $n\gess 4$.
\end{icor}

These do not hold if the assumption~(WH) in Theorem~\ref{T2} is unsatisfied, see \cite{ex}. Note that the condition $\kod\notin\R_Z$ follows from the inequality $k>(n{-}1)d\mi n$ if $n\eq3$, using for instance \cite[Corollary 3.6]{dFEM}.
\sk
Let $\tau_Z:=\sum_{z\in{\rm Sing}\,Z}\tau_{h_z}$ with $\tau_{h_z}$ the {\it Tjurina number} of a local defining function $h_z$, that is, $\tau_{h_z}=\dim\OO_{Y,z}/\bl((\dd h_z),h_z\br)$. Under the assumptions~(IS)-(WH), this coincides with the Milnor number $\mu_{h_z}:=\dim\OO_{Y,z}/(\dd h_z)$. We have the inequalities (see \cite{kosz}):
\begin{equation} \label{8}
\aligned&\mu''_k:=\dim M''_k\les\tau_Z,\q\nu_k:=\dim N_k\les\tau_Z\q(\forall\,k\in\Z),\\&\h{where the equalities hold if}\,\,\,k\ges nd.\endaligned
\end{equation}
\sk
The differential $\ddd$ of $\Om^{\ssb}$ induces $\ddd^{(1)}:N_{k+d}\to M_k$ ($k\in\Z$) so that $M^{(2)}_k$, $N^{(2)}_{k+d}$ are respectively its cokernel and kernel. It is easy to calculate the dimensions of $M_k$, $N_{k+d}$ by using computers, and moreover, under the assumptions~(IS)-(WH), we have the following:
\begin{equation} \label{9}
\h{The composition}\,\,\,N_{k+d}\buildrel{\ddd^{(1)}}\over\longrightarrow M_k\onto M''_k\,\,\,\,\h{is {\it injective} if}\,\,\,\,\kod\notin\R_Z.
\end{equation}
This follows from \cite[Theorem 5.3 and Remark 5.6(i)]{kosz} together with the assertion that the Hodge and pole order filtrations coincide in the case of {\it weighted homogeneous isolated singularities,} see \eqref{7.2} below.

\begin{irem} \label{R2}
In order to determine $\R_f^0$, we do {\it not have to calculate} $\ddd^{(1)}:N_{k+d}\to M_k$ in \eqref{9} for $\kod\in\R_Z$ because of the last inclusion $\R_Z\subset\R_f$ in \eqref{1}. In particular, the information of $H^{n-2}(\Ff,\C)$ (which is given by the kernel of $\ddd^{(1)}$) is unnecessary to determine $\R_f^0$.
\end{irem}

Let $(\dd f)\subset R$ be the {\it Jacobian ideal\,} generated by the partial derivatives of $f$. Then $M$ is identified with the {\it Jacobian ring\,} $R/(\dd f)$ with grading {\it shifted by\,} $-n$.
It is well known that the $\nu_k=\dim N_k$ are expressed by the $\mu_k=\dim M_k$; more precisely
\begin{equation} \label{10}
\nu_k=\mu_k-\ga_k\q(\forall\,k\in\Z)\q\q\h{with}\q\q
\msum_{k\in\Z}\,\ga_{k\,}t^{\,k}=\bl(t+\cdots+t^{\,d-1}\br){}^n,
\end{equation}
see \cite[Formula (3)]{kosz}. (This follows from the assertion that the Euler characteristic of a bounded complex of finite dimensional vector spaces is independent of the differential.)
\sk
By definition we have $\nu_{d+n}\nes0$ if and only if $f$ is annihilated by a nonzero vector field $\xi$ of degree $0$. Here ``degree 0" means that the coefficients of $\xi$ are linear combinations of coordinates of $\C^n$, that is, $\xi\eq\msum_{i,j}\,c_{i,j}x_i\dd_{x_j}$ ($c_{i,j}\ins\C)$, see also Remark~\ref{R8.4} below. The second main theorem of this paper is as follows.

\begin{ithm} \label{T4}
Under the assumptions~{\rm(IS)}-{\rm(WH)}, assume further that $n\eq3$ and $\nu_{d+3}\eq0$, that is, $f$ is not annihilated by any nonzero vector field of degree $0$. Then we have $\mu_k>\nu_{d+k}$, hence $\kod\in\R_f$, for any $k\in\Z\cap[3,d{-}1]$.
\end{ithm}

This follows from Theorem~\ref{T2} together with \cite[Theorem 2]{bcm} (that is, Theorem~\ref{T1} in this paper) and also the non-degeneracy of the Grothendieck residue pairing, see \cite[p.\,659]{GH}. We have {\it symmetries\one} of $\mu'_k:=\dim M'_k$ and $\de''_k:=\mu''_k\mi\nu_{k+d}$ with $\mu''_k:=\dim M''_k$, whose centers are $nd/2$ and $(n{-}1)d/2$ respectively, that is,
\begin{equation} \label{11}
\mu'_k\eq\mu'_{nd-k},\q\de''_k\eq\de''_{(n-1)d-k}\q(\forall\,k\ins\Z),
\end{equation}
see \cite[Corollaries 1 and 2]{kosz} (and also \cite{SeE} for the first equality). Using this and the Lefschetz property of $M'_{\ssb}$ (in particular, $\mu'_{k-1}\less\mu'_k$ for $k\less\tfrac{3d}{2}$ with $n\eq3$, see \cite{DiPo}), Theorem~\ref{T4} and Corollary~\ref{C3} imply the following.

\begin{icor} \label{C4}
Under the hypotheses of Theorem~{\rm\ref{T4}}, we have
\begin{equation} \label{12}
\R_f=\tfrac{1}{d}(\Z\cap[3,k''])\cup\R_Z,
\end{equation}
with $k''=\max(2d{-}3,k'_{\rm max})$, $\,k'_{\rm max}:=\max\{k\mid\mu'_k\nes0\}$.
\end{icor}

To show the inclusion $d(\R_f\stm\R_Z)\sst[3,+\infty)$, we need \cite[Theorem 2.2]{bcm} together with \cite[Corollary 3.6]{dFEM}. In the case of essential indecomposable reduced central hyperplane arrangements in $\C^3$, Corollary~\ref{C4} is shown in \cite{Ba} by a completely different method. Here we have $\max\R_f\slt2\mi\tfrac{1}{d}$, hence $k''\less 2d{-}2$, see \cite[Theorem 1]{bha}. In the case $n\eq3$, the latter theorem is {\it reproved\one} in \cite[Corollary 7.3]{DiSt2} using an estimate of Castelnuovo-Mumford regularity in \cite[Corollary 3.5]{DIM}, see also \cite{Ba}. We can verify that $f$ is not annihilated by any nonzero vector field of degree $0$ for a hyperplane arrangement as above by applying the following (which is a consequence of Propositions~\ref{P8.1} and \ref{P8.2} below).

\begin{ithm} \label{T5}
Under the assumptions~{\rm(IS)}-{\rm(WH)}, assume further $n\eq3$, $\nu_{d+3}\nes0$, and $d\gess5$. Then $f$ is extremely degenerated, see Remark~
{\rm\ref{R8.1}} below.
\end{ithm}

Indeed, assuming that an extremely degenerated $Z$ is essential and indecomposable, that is, $f$ is not equal to $g(x,y)$ or $g(x,y)z$ for any coordinate system $(x,y,z)$, we can verify that the second minimal spectral number at some singular point of $Z$ has multiplicity one at least if $d\gess 4$.
\sk
By Corollary\,\,\ref{C4} and Theorem\,\,\ref{T5}, the roots of Bernstein-Sato polynomial for an essential indecomposable reduced line arrangement $Z$ in $\PP^2$ is determined only by the degree $d\eq\deg Z$ and the local multiplicities ${\rm mult}_pZ$ for $p\ins{\rm Sing}\,Z$ together with vanishing or nonvanishing of the dimension $\mu'_{2d-2}\,({=}\,\mu'_{d+2})$ of the $y$-torsion part $M'_{2d-2}$, where the ``combinatorial data" about the relation between the global components and the local singularities are {\it not\one} needed, see also \cite{Ba}.
Note that the roots of Bernstein-Sato polynomial for a reduced homogeneous polynomial of two variables with degree $d$ are given up to sign by $k/d$ for $k\ins[2,2d{-}2]\cap\Z$. (This follows for instance from \cite{Ma1}.)
\sk
As a conclusion it is recommended to try to run the code of A.\ref{A.3} in Appendix below to calculate the $\mu'_{k/d}$ and $\de_{k/d}\eq\mu'_{k/d}\pl\de''_{k/d}$ in the case conditions~(IS)-(WH) are satisfied with $n\eq3$. If the coefficients of $T^k$ in the last output about $\de_{k/d}$ are all positive for $k\in[3,d{-}1]$, then $d\one\R_f\cap\Z$ coincides with $\{3,\dots,k''\}$ outside $d\one\R_Z$ where $k''$ is as in Corollary~\ref{C4}. (Note that a polynomial with $\nu_{d+3}\nes0$ can produce lots of examples of the same type applying coordinate changes.) The coefficient of $T^3$ of the last output is less than or equal to $0$ only in the case $\nu_{d+3}\nes0$. There are essentially only two examples of $\nu_{d+3}$-nonvanishing reduced plane curves (up to deformation) which have only weighted homogeneous singularities and are {\it not\one} extremely degenerated, where $d\eq3,4$; the singularities become non-weighted-homogeneous if $d\gess5$, see Propositions~\ref{P8.1}--\ref{P8.2} and Remark~\ref{R8.5} below.
\sk
In Part~\ref{Pa1} we review some basics of pole order spectral sequences, and give the proofs of Proposition~\ref{P1} and Theorem~\ref{T5.1} on the compatibility of certain self-duality isomorphisms. In Part~\ref{Pa2} we study the fundamental exact sequences about the vanishing cycles, and prove Theorem~\ref{T10.1} which is a key to the proof of Theorem~\ref{T2}. In Part~\ref{Pa3} we calculate the filtered twisted de Rham complexes of weighted homogeneous polynomials with isolated singularities. In Part~\ref{Pa4} we prove Theorems~\ref{T2}--\ref{T4} after recalling some basics of Bernstein-Sato polynomials. In Part~\ref{Pa5} we calculate some examples explicitly. In Appendix we explain how to use the computer programs Macaulay2 and Singular for explicit computations of roots of Bernstein-Sato polynomials.
\sk
We thank A.~Dimca for useful discussions on pole order spectral sequences and for making a computer program based on \cite{nwh} (which encouraged us to make a less sophisticated one).
This work was partially supported by Kakenhi 24540039 and 15K04816.

\tableofcontents
\numberwithin{equation}{section}

\part{Preliminaries} \label{Pa1}
\nin
In this part we review some basics of pole order spectral sequences, and give the proofs of Proposition~\ref{P1} and Theorem~\ref{T5.1} on the compatibility of certain self-duality isomorphisms.

\section{Algebraic microlocal Gauss-Manin systems} \label{S1}
For a homogeneous polynomial $f$, we have the {\it algebraic microlocal Gauss-Manin complex} $\Ct^{\ssb}_f$ defined by
\begin{equation} \label{1.1}
\Ct_f^j:=\Om^j[\dd_t,\dti]\q\h{with differential}\q\ddd-\dd_t\,\dfw,
\end{equation}
where $\Om^{\ssb}$ is as in the \nameref{intr}.
It is a complex of graded $\C[t]\langle\dd_t,\dti\rangle$-modules with
\begin{equation} \label{1.2}
\deg t=-\deg\dd_t=d,
\end{equation}
and the actions of $t$, $\dd_t^i$ are defined by
\begin{equation} \label{1.3}
t(\om\dd_t^k)=(f\om)\dd_t^k-k\,\om\dd_t^{k-1},\q\dd_t^i (\om\dd_t^k)=\om\dd_t^{k+i}\q(\om\in\Om^j,\,i\in\Z).
\end{equation}
\sk
The {\it algebraic microlocal Gauss-Manin systems} are defined by
$$\Gt_f^j:=H^{n-j}\Ct^{\ssb}_f\q\q(j\in[0,\sigma_Z+1]),$$
where $\sigma_Z:=\dim{\rm Sing}\,Z$ as in the \nameref{intr}. These are free graded $\C[\dd_t,\dti]$-modules of finite type, and there are canonical isomorphisms
\begin{equation} \label{1.4}
\Gt^j_{f,k}=\Ht^{n-1-j}(\Ff,\C)_{\ee(-k/d)}\q\q\bl(\forall\,k\in\Z\br),
\end{equation}
where $\ee(\al):=\exp(2\pi i\al)$ for $\al\in\Q$, $\Ff:=f^{-1}(1)$ (which is viewed as the Milnor fiber of $f$), and $E_{\la}:={\rm Ker}(T_s-\la)$ in $E:=\Ht^{\ssb}(\Ff,\C)$ with $T_s$ the semisimple part of the monodromy $T$, see \cite{ggr}, \cite{kosz}, etc.

\section{Pole order spectral sequences} \label{S2}
In the notation of Section~\ref{S1}, we have the filtration $P'$ on $\Ct_f^{\ssb}$ defined by
\begin{equation} \label{2.1}
P'_i\,\Ct_f^j:=\mopl_{k\les i+j}\,\Om^j\dd_t^k.
\end{equation}
This is an exhaustive increasing filtration, and induces the filtration $P'$ on the Gauss-Manin systems $\G_f^j$ compatible with the grading. Since $\deg\dd_t=-d$, we have
\begin{equation} \label{2.2}
\Om^j\dd_t^k=\Om^j(kd),
\end{equation}
where $(m)$ for $m\in\Z$ denotes the shift of grading as in Proposition~\ref{P1}. Set
$$P^{\prime\,i}=P'_{-i}.$$
We have
\begin{equation} \label{2.3}
\Gr^i_{P'}\Ct_f^{\ssb}=K_f^{\ssb}((n{-}i)d).
\end{equation}
Here $K_f^{\ssb}$ is the {\it Koszul complex} defined by
\begin{equation} \label{2.4}
K_f^j:=\Om^j((j{-}n)d)\q\h{with differential}\q\dfw.
\end{equation}
We have the {\it microlocal pole order spectral sequence}
\begin{equation} \label{2.5}
{}_{P'\!}\Et_1^{i,j-i}=H^j\Gr^i_{P'}\Ct_f^{\ssb}\Longrightarrow\Gt_f^{j+n}\,\bl(=H^j\Ct_f^{\ssb}\br).
\end{equation}
In the notation of Proposition~\ref{P1}, we have
\begin{equation} \label{2.6}
{}_{P'\!}\Et_r^{i,j-i}=\begin{cases}\bl(H^j_{\dfw}\Om^{\ssb}\br)((j{-}i)d)=\bl(H^jK_f^{\ssb}\br)((n{-}i)d)&\h{if}\,\,\,r=1,\\ H^j_{\ddd}\bl(H^{\ssb}_{\dfw}\Om^{\ssb}\br)((j{-}i)d)&\h{if}\,\,\,r=2.\end{cases}
\end{equation}
Set
\begin{equation} \label{2.7}
P:=P'[1],\q\h{that is,}\q P^i:=P^{\prime\,i+1}\q(i\in\Z).
\end{equation}
The isomorphisms in \eqref{1.4} induce the filtered isomorphisms
\begin{equation} \label{2.8}
\bl(\Gt^j_{f,k},P\br)=\bl(\Ht^{n-1-j}(\Ff,\C)_{\ee(-k/d)},P\br)\q\q\bl(\forall\,k\in[1,d]\br),
\end{equation}
where $P$ on the right-hand side is the {\it pole order filtration}, see \cite{Di1}, \cite[Section 1.8]{ggr}.
\sk
By the definition \eqref{2.1} we have the filtered graded isomorphisms
\begin{equation} \label{2.9}
\dd_t^k:\bl(\Ct_f^{\ssb},P'\br)\simto\bl(\Ct_f^{\ssb}(-kd),P'[-k]\br)\q(k\in\Z).
\end{equation}
Recall that $(P'[m])^i=P'{}^{i+m}$, $(P'[m])_i=P'_{i-m}$, and $G(m)_k=G_{m+k}$ for any graded module $G_{\ssb}$ and $m\in\Z$ in general.
\sk
By \eqref{2.9} we get the graded isomorphisms of spectral sequences for $k\in\Z:$
\begin{equation} \label{2.10}
\dd_t^k:{}_{P'\!}\Et_r^{i,j-i}\simto{}_{P'\!}\Et_r^{i-k,j-i+k}(-kd)\q(r\ges 1),
\end{equation}
which are compatible with the differentials $\ddd_r$ of the spectral sequences.

\section{Isolated singularity case} \label{S3}
In the notation of the \nameref{intr}, assume $Z\subset Y$ has only isolated singularities, that is,
\begin{equation} \label{3.1}
\sigma_Z=\dim\Si=0\q\h{with}\q\Si={\rm Sing}\,Z.
\end{equation}
Set
\begin{equation} \label{3.2}
M:=H^nK_f^{\ssb},\q N:=H^{n-1}K_f^{\ssb},
\end{equation}
where $k_f^{\ssb}$ is as in \eqref{2.4}. These are graded $\C$-vector spaces.
\sk
Let $y$ be a sufficiently general linear combination of coordinates $x_1,\dots,x_n$ of $\C^n$. Then $M$, $N$ are finitely generated graded $\C[y]$-modules.
It is known that $N$ is $y$-torsion-free (see, for instance, \cite{kosz}). Set
\begin{equation} \label{3.3}
M':=M_{\rm tor},\q M'':=M_{\rm free}=M/M_{\rm tor},
\end{equation}
where $M_{\rm tor}$ denotes the $y$-torsion part of $M$. (The latter is independent of $y$ as long as $y$ is sufficiently general). These are also finitely generated graded $\C[y]$-modules. Set
\begin{equation} \label{3.4}
\aligned\Mt&:=M[y^{-1}]=M\otimes_{\C[y]}\C[y,y^{-1}]\,\bl(=\Mt'':=M''[y^{-1}]\br),\\\Nt&:=N[y^{-1}]=N\otimes_{\C[y]}\C[y,y^{-1}].\endaligned
\end{equation}
These are finitely generated free graded $\C[y,y^{-1}]$-modules.
By the grading, we have the direct sum decompositions
$$M=\mopl_{k\in\Z}\,M_k,\q N=\mopl_{k\in\Z}\,N_k,\q\Mt=\mopl_{k\in\Z}\,\Mt_k,\q\Nt=\mopl_{k\in\Z}\,\Nt_k,\,\,\h{etc.}$$
By definition there are isomorphisms
\begin{equation} \label{3.5}
y:\Nt_k\simto\Nt_{k+1},\q y:\Mt_k\simto\Mt_{k+1}\q(k\in\Z),
\end{equation}
together with natural inclusions
\begin{equation} \label{3.6}
N_k\into\Nt_k,\q M''_k\into\Mt''_k=\Mt_k\q(k\in\Z),
\end{equation}
inducing isomorphisms for $k\gg 0$, since $N$ and $M''$ are $y$-torsion-free.
In the notation of \cite[5.1]{kosz} we have
\begin{equation} \label{3.7}
H^{-1}(\spKf)=\Nt,\q H^0(\spKf)=\Mt.
\end{equation}
\sk
In the notation of Section~\ref{S2}, set
\begin{equation} \label{3.8}
M^{(r)}:={}_{P'\!}\Et_r^{n,0},\q N^{(r)}:={}_{P'\!}\Et_r^{n,-1}\q(r\ges 1).
\end{equation}
Using the differentials $\ddd_r$ of the microlocal pole order spectral sequence \eqref{2.5} together with the isomorphisms in \eqref{2.10}, we get the morphisms of graded $\C$-vector spaces of degree $-rd:$
\begin{equation} \label{3.9}
\ddd^{(r)}:N^{(r)}\to M^{(r)}\q\h{for any}\,\,\,r\ges 1,
\end{equation}
such that $N^{(r)}$, $M^{(r)}$ are respectively the kernel and cokernel of $\ddd^{(r-1)}$ for any $r\ges 2$, and are independent of $r\gg 0$ (that is, $\ddd^{(r)}=0$ for $r\gg 0)$. More precisely, \eqref{3.9} is given by the composition (using the isomorphisms \eqref{2.10}):
\begin{equation} \label{3.10}
{}_{P'\!}\Et_r^{n,0}\buildrel{\ddd_r}\over\longrightarrow{}_{P'\!}\Et_r^{n+r,-r-1}\buildrel{\dd_t^r}\over{\longrightarrow}{}_{P'\!}\Et_r^{n,-1}(-rd).
\end{equation}
\sk
By \eqref{2.3} or \eqref{2.6}, we have for $r=1$
\begin{equation} \label{3.11}
M^{(1)}=M,\q N^{(1)}=N,
\end{equation}
and $\ddd^{(1)}:N\to M$ is induced by the differential $\ddd$ of the de Rham complex $(\Om^{\ssb},\ddd)$. By \eqref{2.6} we get moreover
\begin{equation} \label{3.12}
M^{(2)}=H^n_{\ddd}\bl(H^{\ssb}_{\dfw}\Om^{\ssb}\br),\q N^{(2)}=H^{n-1}_{\ddd}\bl(H^{\ssb}_{\dfw}\Om^{\ssb}\br)(-d).
\end{equation}
These $M^{(r)}$, $N^{(r)}$, and $\ddd^{(r)}$ are equivalent to the microlocal pole order spectral sequence because of the isomorphisms \eqref{2.9}--\eqref{2.10}. Set
\begin{equation} \label{3.13}
M^{(\infty)}:=M^{(r)},\q N^{(\infty)}:=N^{(r)}\q(r\gg 0).
\end{equation}
Here $\mu^{(r)}_k:=\dim M^{(r)}_k$, $\nu^{(r)}_k:=\dim N^{(r)}_k$ are finite and non-increasing for $r\ges 1$ with $k$ fixed. Hence they are stationary for $r\gg 0$ with $k$ fixed.
\sk
Note that the shift of the grading by $-d$ for $N^{(2)}$ in \eqref{3.12} comes from the definition of the Koszul complex $K_f^{\ssb}$ where the differential $\dfw$ preserves the grading.

\section{Weighted homogeneous isolated singularity case} \label{S4}
Assume the isolated singularity condition (IS) together with (WH) in Theorem~\ref{T2}, see also Remark~\ref{R4.1} below. In the notation of Section~\ref{S3}, set as in \cite[5.1]{kosz}:
\begin{equation} \label{4.1}
h:=f/y^d\q\h{on}\q Y':=Y\setminus\{y=0\}\cong\C^{n-1}.
\end{equation}
Note that $h$ is a defining function of $Z\setminus\{y=0\}$ in $Y'$, and we have $\Si={\rm Sing}\,Z\subset Y'$ (since $y$ is sufficiently general).
\sk
For $z\in\Si$, set
\begin{equation} \label{4.2}
\Xi_{h_z}:=\Om_{Y'_{\rm an},z}^{n'}/\ddd h_z{\wedge}\Om_{Y'_{\rm an},z}^{n'-1}\q(n':=\dim Y'=n{-}1).
\end{equation}
These are finite dimensional vector spaces of dimension $\mu_z$, where $\mu_z$ is the Milnor number of $(h_z,z)$ which coincides with the Tjurina number $\tau_z$ ($z\in\Si$) under the assumption~(WH) in Theorem~\ref{T2}.
We have the canonical isomorphism
\begin{equation} \label{4.3}
\Xi_{h_z}^{\rm alg}\simto\Xi_{h_z},
\end{equation}
where $\Xi_{h_z}^{\rm alg}$ is defined in the same way as $\Xi_{h_z}$ with $\Om_{Y'_{\rm an},z}^{\ssb}$ replaced by $\Om_{Y',z}^{\ssb}$. (Indeed, $\OO_{Y'_{\rm an},z}$ is flat over $\OO_{Y',z}$ so that the tensor product of $\OO_{Y'_{\rm an},z}$ over $\OO_{Y',z}$ is an exact functor. Applying this to the exact sequence defining $\Xi_{h_z}^{\rm alg}$, that is,
$$\Om_{Y',z}^{n'-1}\buildrel{\ddd h{\wedge}}\over\longrightarrow\Om_{Y',z}^{n'}\to\Xi_{h_z}^{\rm alg}\to 0,$$
we get the isomorphism
$$\Xi_{h_z}=\Xi_{h_z}^{\rm alg}\,{\otimes}_{\OO_{Y',z}}\OO_{Y'_{\rm an},z}.$$
We have a canonical morphism $\Xi_{h_z}^{\rm alg}\to\Xi_{h_z}$ induced by the inclusion $\OO_{Y',z}\into\OO_{Y'_{\rm an},z}$. So the assertion follows by suing a finite filtration on $\Xi_{h_z}^{\rm alg}$ such that its graded quotients are $\C$, since $\C\otimes_{\OO_{Y',z}}\OO_{Y'_{\rm an},z}=\C$.)
\sk
Fixing $y$, we have the isomorphisms compatible with \eqref{3.5}
\begin{equation} \label{4.4}
\Nt_k=\mopl_{z\in\Si}\,\Xi_{h_z},\q\Mt_k=\mopl_{z\in\Si}\,\Xi_{h_z}\q(\forall\,k\in\Z).
\end{equation}
This follows from the argument in \cite[5.1]{kosz} by using the isomorphisms in \eqref{3.7}, see also the proof of Theorem~\ref{T5.1} below.
\sk
By the theory of Gauss-Manin connections on Brieskorn lattices (see \cite{Br}) in the weighted homogeneous polynomial case, we have the finite direct sum decompositions
\begin{equation} \label{4.5}
\Xi_{h_z}=\mopl_{\al\in\Q}\,\Xi_{h_z}^{\al}\q(z\in\Si).
\end{equation}
defined by
\begin{equation} \label{4.6}
\Xi_{h_z}^{\al}:={\rm Ker}(\dd_tt-\al)\subset\Xi_{h_z}\q(\al\in\Q,\,z\in\Si),
\end{equation}
which is independent of the choice of analytic local coordinates. (Note that we have a well-defined action of $\dd_tt$ on $\Xi_{h_z}$, since $h_z$ is contained in the ideal generated by its partial derivatives, see \eqref{17.2} below.)

\begin{rem} \label{R4.1}
We say that $h$ is a {\it weighted homogeneous polynomial} with positive weights $(w_i)$ for a coordinate system $(y_i)$ of $Y:=\C^n$, if $h$ is a linear combination of $\mprod_iy_i^{a_i}$ with $\sum_iw_ia_i=1$, where $w_i\in\Q_{>0}$, $a_i\in\N$. We have
\begin{equation} \label{4.7}
\la^m h(y_1,\dots,y_n)=h(\la^{m_1}y_1,\dots,\la^{m_n}y_n)\q(\la\in\C^*),
\end{equation}
where $m_i:=mw_i\in\N$ with $m$ the smallest positive integer satisfying $mw_i\in\N$, and \eqref{4.7} holds with $\la$ replaced by $g\in\OO_{Y,0}$. This implies that condition~(WH) in Theorem~\ref{T2} is independent of the choice of a defining function $h_z$.
\sk
We also have
\begin{equation} \label{4.8}
h=v(h)\q\h{with}\q v:=\msum_i\,w_i\,y_i\dd_{y_i}\in\Theta_Y.
\end{equation}
This implies that $h$ is {\it quasihomogeneous}, that is, $h\in(\dd h)$ with the notation in a remark after Theorem~\ref{T2}. The converse holds in the isolated singularity case, see \cite{SaK}.
\end{rem}

\section{Compatibility of duality isomorphisms} \label{S5}
In this section we show the following.

\begin{thm} \label{T5.1}
Under the isomorphisms in {\rm\eqref{4.4}}, the duality isomorphism between $\Nt_k$ and $\Mt_{nd-k}$ induced from the self-duality isomorphism for the Koszul complex in \cite[Thm.\,1]{kosz} $($using the graded local duality as in \cite[1.1.4 and 1.7.3]{kosz}$)$ is identified up to a constant multiple with the direct sum of the canonical self-duality isomorphisms of $\Xi_{h_z}$ for $z\in\Si$.
\end{thm}

\begin{proof}
The self-duality isomorphism used in the proof of \cite[Theorem 1]{kosz} is a canonical one, and is induced by the canonical graded isomorphism
$$\Om^j={\rm Hom}_R(\Om^{n-j},\Om^n),$$
where $\Om^n[n]$ is the graded dualizing complex, see also \cite{Ei}, \cite{Ha1}, etc. (There is no shift of grading here, and it appears in the duality isomorphisms in \cite[Theorem 1]{kosz}.)
\sk
The duality isomorphisms are compatible with the localization by $y$ in \eqref{3.4} (see also \cite[5.1]{kosz}). Here we consider the graded modules only over $\C[y]$ or $\C[y,y^{-1}]$ as in \cite[Remark 1.7]{kosz}, and then consider the graded duals also over $\C[y]$ or $\C[y,y^{-1}]$.
The graded dual over $\C[y]$ is defined by using the graded dualizing complex $(\C[y]\ddd y)[1]$, where the degree of complex is shifted by 1, and the degree of grading is also shifted by 1 because of $\ddd y$.
We have the same with $\C[y]$ replaced by $\C[y,y^{-1}]$.
\sk
As for the relation with the graded local duality as in \cite[1.1.4 and 1.7.3]{kosz}, we can take graded free generators $\{u_i\}$, $\{v_i\}$ of $M''$, $N$ over $\C[y]$ such that $u_i\in M''_{k_i}$, $v_i\in N_{nd+1-k_i}$, with $k_i\in\Z_{>0}$, and $\C[y]u_i$ is orthogonal to $\C[y]v_j$ for $i\ne j$. (Here the shift of grading by 1 comes from the degree of $\ddd y$ in the above remark.) We can get the $u_i$ by using the filtration on $M_k$ defined by $M_k\cap y^pM''$ ($p\in\Z$), and similarly for the $v_i$. Then $\{y^{p-k_i}u_i\}$ and $\{y^{k_i-p-nd-1}u_i\}$ for $p\gg 0$ are $\C$-bases of $\Mt_p=M''_p$ and $\Nt_{-p}=(\Nt/N)_{-p}$, which are orthogonal to each other under the pairing given by \cite[1.1.4 and 1.7.3]{kosz}.
\sk
By the above argument, it is enough to consider the graded $\C[y,y^{-1}]$-dual of $\Mt$ instead of the $\C$-dual of the $M''_p$ for $p\gg 0$.
In particular, we can neglect the shift of grading in the self-duality isomorphism (by using the isomorphisms in \eqref{3.5}) for the proof of Theorem~\ref{T5.1}.
\sk
Consider the blow-up along the origin
$$\pi:\Xt\to X:=\C^n.$$
(This blow-up is used only for the {\it coordinate change} as is explained below. Here we may assume $y=x_n$ replacing $x_n$ if necessary.)
\sk
Let $\Xt''\subset\Xt$ be the complement of the {\it total} transform of $\{y=0\}\subset X$. We have the isomorphism
\begin{equation} \label{5.1}
\Xt''=\C^{n'}\times\C^*\q\h{with}\q n'=n{-}1.
\end{equation}
(This must be distinguished from the isomorphism $\Xt''=X\setminus\{y=0\}\cong\C^{n'}\times\C^*$, where the effect of the blow-up is neglected, since the coordinates $x_1,\dots,x_n$ are used here instead of $x'_1,\dots,x'_{n'}$ defined below.)
\sk
In the notation of Section~\ref{S3}, set $x'_i:=x_i/y$ ($i\in[1,n']$), and
$$R':=\C[x'_1,\dots,x'_{n'}].$$
Then
\begin{equation} \label{5.2}
\pi^*f|_{\Xt''}=y^dh\q\h{with}\q h:=f/y^d\in R',
\end{equation}
see also \eqref{4.1}. This decomposition is compatible with the decomposition \eqref{5.1}, where $x'_1,\dots,x'_{n'}$ and $y$ are identified with the coordinates of $\C^{n'}$ and $\C^*$ respectively. Since
$$\dd_y(y^dh)=dy^{d-1}h,\q\dd_{x'_i}(y^dh)=y^dh_i\q\h{with}\q h_i:=\dd_{x'_i}h\in R',$$
the localized Koszul complex $\Kt_{y^dh}^{\ssb}$ for the partial derivatives of $y^dh$ on $\Xt''=Y'\times\C^*$ can be identified, up to a non-zero constant multiple, with the shifted mapping cone
\begin{equation} \label{5.3}
C\bl(h:K_h^{\ssb}[y,y^{-1}]\to K_h^{\ssb}[y,y^{-1}]\br)[-1],
\end{equation}
where $K_h^{\ssb}[y,y^{-1}]$ is the scalar extension by $\C\into\C[y,y^{-1}]$ of the Koszul complex $K_h^{\ssb}$ defined by
$$h_i\in{\rm End}(R')\,\,\,\,(i\in[1,n']).$$
Note that $y^d\ddd y$ and $y^d$ can be omitted in the above descriptions, since they express only the shift of grading, and can be neglected by using the isomorphisms in \eqref{3.5} as is explained above. (Indeed, only $y$ has non-zero degree, since $\deg y=1$ and $\deg h=\deg h_i=\deg x'_i=0$.)
\sk
By \eqref{4.3} there is a quasi-isomorphism (or an isomorphism in the derived category)
\begin{equation} \label{5.4}
K_h^{\ssb}[y,y^{-1}]\simto\mopl_{z\in\Si}\,\Xi_{h_z}[y,y^{-1}][-n'].
\end{equation}
\sk
Consider the filtration $G$ on $\Kt_{y^dh}^{\ssb}$ such that
\begin{equation} \label{5.5}
\aligned\Gr_G^0\Kt_{y^dh}^{\ssb}&\cong K_h^{\ssb}[y,y^{-1}]\,\bl(\simto\mopl_{z\in\Si}\,\Xi_{h_z}[y,y^{-1}][-n']\br),\\ \Gr_G^1\Kt_{y^dh}^{\ssb}&\cong K_h^{\ssb}[y,y^{-1}][-1]\,\bl(\simto\mopl_{z\in\Si}\,\Xi_{h_z}[y,y^{-1}][-n'{-}1]\br),\endaligned
\end{equation}
where $\Kt_{y^dh}^{\ssb}$ is identified with the mapping cone \eqref{5.3} (and \eqref{5.4} is used for the second isomorphisms in the derived category).
\sk
By condition (WH) in Theorem~\ref{T2}, we have the vanishing of the morphisms
$$h:\Xi_{h_z}[y,y^{-1}]\to\Xi_{h_z}[y,y^{-1}]\q(\forall\,z\in\Si).$$
So the filtration $G$ on $\Kt_{y^dh}^{\ssb}$ {\it splits} in the bounded derived category of graded $R'[y,y^{-1}]$-modules $D^b(R'[y,y^{-1}])_{\rm gr}$, and there is a (non-canonical) isomorphism
\begin{equation} \label{5.6}
\Kt_{y^dh}^{\ssb}\cong\Gr_G^0\Kt_{y^dh}^{\ssb}\oplus\Gr_G^1\Kt_{y^dh}^{\ssb}\q\h{in}\,\,\,D^b(R'[y,y^{-1}])_{\rm gr},
\end{equation}
compatible with the filtration $G$ so that the $\Gr_G^k$ of \eqref{5.6} are the identity morphisms for $k\eq0,1$. Here we use the isomorphism
$${\rm Hom}_{\mathcal A}(M,M')={\rm Hom}_{D^b({\mathcal A})}(M,M'),$$
for objects $M,M'$ of an abelian category ${\mathcal A}$ in general.
\sk
The self-duality isomorphism for $\Kt_{y^dh}^{\ssb}$ is compatible with the filtration $G$, and induces a duality between
$$\Gr_G^0\Kt_{y^dh}^{\ssb}\q\h{and}\q\Gr_G^1\Kt_{y^dh}^{\ssb},$$
which can be identified, up to a non-zero constant multiple and also a shift of grading, with the canonical self-duality of $K_h^{\ssb}[y,y^{-1}]$ by the above argument.
Moreover the last self-duality is the {\it scalar extension} by $\C\into\C[y,y^{-1}]$ of the self-duality of the Koszul complex $K_h^{\ssb}$.
So the assertion follows (since the above induced duality between $\Gr_G^0\Kt_{y^dh}^{\ssb}$ and $\Gr_G^1\Kt_{y^dh}^{\ssb}$ is sufficient for the proof of Theorem~\ref{T5.1}). This finishes the proof of Theorem~\ref{T5.1}.
\end{proof}

\begin{rem} \label{R5.1}
The self-duality of $\Xi_{h_z}$ is given by using the so-called residue pairing (see for instance \cite{Ha1}) which is compatible with the direct sum decompositions \eqref{4.5}, and implies the duality between $\Xi_{h_z}^{\al}$ and $\Xi_{h_z}^{n'-\al}$.
\end{rem}

\begin{rem} \label{R5.2}
In the notation of the \nameref{intr}, $N_k$ and $M''_{dn-k}$ are orthogonal subspaces to each other by the duality in Theorem~\ref{T5.1}. This is compatible with Corollary~\ref{C3}.
\end{rem}

\section{Proof of Proposition~\ref{P1}} \label{S6}
The assertion follows from the filtered isomorphisms \eqref{2.8} together with the identification of the pole order spectral sequence with the $M^{(r)}$, $N^{(r)}$, $\ddd^{(r)}$ in Section~\ref{S3}. This finishes the proof of Proposition~\ref{P1}.

\section{Steenbrink spectrum and Bernstein-Sato polynomials} \label{S7}
Assume $h$ is a weighted homogeneous polynomial of $n$ variables with weights $w_1,\dots,w_n$ as above. The {\it Steenbrink spectrum} $\Sp(h)=\msum_{i=1}^{\mu_h}\,t^{\al_{h,i}}$ and the {\it spectral numbers} $\al_{h,i}\in\Q$ are defined by
\begin{equation} \label{7.1}
\aligned0<\al_{h,1}\les\cdots\les\al_{h,\mu_h}<n,&\\ \#\bl\{i\,\big|\,\al_{h,i}=\al\br\}=\dim\Gr_F^pH^{n-1}(F_h,\C)_{\la}&\,\bl(=\dim\Xi_h^{\al}\br),\endaligned
\end{equation}
where $p:=[n-\al]$, $\la:=\exp(-2\pi i\al)$, and $\mu_h$ is the Milnor number of $h$, see \cite{St2}, \cite{ScSt}. (Here $\Xi_h^{\al}$ is as in Section~\ref{S4}.) We have moreover
\begin{equation} \label{7.2}
b_h(s)=(s+1)\,\bl[\mprod_{i=1}^{\mu_h}(s+\al_{h,i})\br]_{\rm red},
\end{equation}
where $\bl[\mprod_j(s+\beta_j)^{m_j}\br]_{\rm red}:=\mprod_j(s+\beta_j)$ for $\beta_i\ne\beta_j$ ($i\ne j$) and $m_j\in\Z_{>0}$ in general. This equality can be proved by combining \cite{Ma1} and \cite{ScSt}, \cite{Va1} (which show that the Bernstein-Sato polynomial and the Hodge filtration on the Milnor cohomology can be obtained by using the Brieskorn lattice \cite{Br}), see also \cite{Sat}, \cite{St1}. It is also well known that
\begin{equation} \label{7.3}
\Sp(h)=\mprod_{j=1}^n\,(t^{w_j}-t)/(1-t^{w_j}),
\end{equation}
see for instance \cite{JKSY0}. We may use this to deduce the {\it symmetry of spectral numbers}
\begin{equation} \label{7.4}
\al_{h,i}=\al_{h,j}\q\q(i+j=\mu_h+1).
\end{equation}
Taking the limit of \eqref{7.3} for $t\to 1$, we can get
\begin{equation} \label{7.5}
\mu_h=\mprod_{j=1}^n\,\bl(\tfrac{1}{w_j}-1\br).
\end{equation}
This is well known in the Brieskorn-Pham type case, that is, if $w_i=a_i^{-1}$ with $a_i\in\N$. For any homogeneous polynomial of degree $m$ having an isolated singularity at 0, the assertions \eqref{7.2} and \eqref{7.3} imply that the roots of $b_f(s)$ are $\tfrac{i}{m}$ up to sign for $i\ins[n,nm{-}n]\cap\Z$, where $a_i\eq m$ ($\forall\,i$).

\begin{rem} \label{R7.1}
We can prove \eqref{7.3} by using the Koszul complex $K_h^{\ssb}$ as in \eqref{2.4} with $f$ replaced by $h$. Indeed, the denominator of the right hand side of \eqref{7.3} gives the Hilbert series of the graded ring $R$ with $\deg x_i=w_i$, and the numerator corresponds to the shift of the grading of each component of the Koszul complex (since $\deg \dd_{x_i}h=1-w_i$), where the complex is viewed as the associated single complex of an $n$-ple complex associated with the multiplications by partial derivatives $\dd_{x_i}h$ ($i\in[1,n]$), see for instance \cite[Section IV.2]{Se} (and also \cite{JKSY0}).
\sk
We can use \eqref{7.3} for an explicit computation using Macaulay2, see A.\ref{A.1} in Appendix below.
\end{rem}

\section{Special hypersurfaces} \label{S8}
A reduced homogeneous polynomial $f$ is called ``special", if for a general $\C$-linear combination $y$ of coordinates $x_1,\dots,x_n$ of $X\defs\C^n$, there is a nonzero $\C$-linear combination $z$ of coordinates such that $zh\in(\dd g)\sst R'$. Here $R'\defs\C[x']$ with $x'\eq(x'_1,\dots,x'_{n-1})$ a coordinate system of $X'_y\defs\{y\eq0\}\sst X$, and $g,h$ are respectively the restrictions of $f$ and $\eta f$ to the hypersurface $X'_y$ with $\eta$ any constant vector field (that is, a $\C$-linear combination of the $\dd_{x_i}$) such that $\eta y\ne 0$.
\sk
For $j\eq0,1$, let $G_{\ssb}^j$ be the $j$\one th cohomology of the shifted mapping cone of the multiplication by $h$\,:
\begin{equation} \label{8.1}
m_h:(R'/(\dd g))(1{-}d)\to R'/(\dd g).
\end{equation}
From the short exact sequence of complexes
\begin{equation} \label{8.2}
0\to K^{\ssb}_f(-1)\buildrel{y}\over\to K^{\ssb}_f\to K^{\ssb}_f/yK^{\ssb}_f(-1)\to 0,
\end{equation}
we get the long exact sequence
\begin{equation} \label{8.3}
0\to N_{\ssb}(n{-}1)\buildrel{y}\over\to N_{\ssb}(n)\to G_{\ssb}^0\to M_{\ssb}(n{-}1)\buildrel{y}\over\to M_{\ssb}(n)\to G_{\ssb}^1\to0,
\end{equation}
using the $n$-ple complex structure of the Koszul complex. Here $M_{\ssb}\eq H^nK_f^{\ssb}$, $N_{\ssb}\eq H^{n-1}K_f^{\ssb}$, and $K^n_f$, $K^{n-1}_f$ are direct sums of copies of $R\defs\C[x]$ shifted respectively by $-n$ and $1{-}n{-}d$, see \eqref{2.4}. These shifts are compatible with the shift by the differential of the complex defining the $G_{\ssb}^j$ which is defined by multiplication by $h$ of degree $d{-}1$, see \eqref{8.1}.
\sk
If $f$ is not ``special" in the above sense, we see that
\begin{equation} \label{8.4}
G^0_k=0\q\h{for}\,\,\,k\less d.
\end{equation}
Indeed, for $z\ins R'_1$, we have $zh\ins(\dd g)$ if and only if $z\ins G^0_d$. By the exact sequence \eqref{8.3}, we then get that
\begin{equation} \label{8.5}
\nu_{d+n}\nes0\,\Longrightarrow\h{$f$ is ``special".}
\end{equation}
Moreover
\begin{equation} \label{8.6}
\h{the converse of \eqref{8.5} holds if $n\eq3$ and $d\gess 5$}.
\end{equation}
Indeed, the Lefschetz property of $M'_{\ssb}$ for $n\eq3$ (see \cite{DiPo}) implies that
\begin{equation} \label{8.7}
{\rm Ker}(y:M_k\to M_{k+1})=0\q\h{if}\,\,\,k\slt\tfrac{3d}{2}.
\end{equation}
Note that $M'(2)_{d}\eq M'_{d+2}$, and $d{+}2\slt\tfrac{3d}{2}\iff d\sgt 4$.

\begin{rem} \label{R8.1}
We say that a reduced hypersurface $Z\sst\PP^{n-1}$ or its defining polynomial $f$ is {\it extremely degenerated\,} if there is a coordinate system $(x_1,\dots,x_n)$ of $\C^n$ such that $f$ is a linear combination of monomials $x^{\nu}$ with $\nu=(\nu_1,\dots,\nu_n)\ins\N^n$ satisfying $\msum_{i=1}^n\,c_i\nu_i\eq0$ for some fixed $(c_1,\dots,c_n)\ins\Q^n$ with at least one of the $c_i$ nonzero. Here the last condition is equivalent to that $\msum_{i=1}^n\,\eta_i\dd_{x_i}f\eq0$ with $\eta_i\eq c_ix_i$. So $\nu_{d+n}\nes0$ for extremely degenerated polynomials.
\end{rem}

\begin{rem} \label{R8.2}
In the case $n\eq3$ with variables $x,y,z$, the support
$${\rm supp}\,f\defs\bl\{(i,j,k)\ins\N^3\mid a_{i,j.k}\nes0\br\}$$
of an extremely degenerated polynomial $f\eq\msum_{(i,j,k)\in\N^3}\,a_{i,j,k}x^iy^jz^k$ must contain two points $(i',j',k')$, $(i'',j'',k'')\ins\N^3$ with $i',j',k''\ins\{0,1\}$ up to permutation of coordinates. Indeed, the support must contain certain points of $\N^3$ corresponding to the {\it three\one} conditions that $f$ is not divisible by any of the three monomials $x^2$, $y^2$, $z^2$, although it must be contained in a line by hypothesis. In the case $f=x^ay^{d-a}\pl z^d$ ($0\slt a\slt d$) we can calculate the reduced Bernstein-Sato polynomial $b_f(s)/(s{+}1)$ using the Thom-Sebastiani type theorem (see \cite{mic}), where the $c_i$ are $a{-}d,a,0$.
\end{rem}

\begin{rem} \label{R8.3}
In the extremely degenerated case with $n\eq3$, any irreducible component of $Z$ is a rational curve, since it is the closure of an irreducible component of $\{u^a\eq v^b\}\sst\C^2\,({\subset}\,\PP^2)$ if it is not a line. The singular points of $Z$ are contained in those of $\{xyz\eq0\}\sst\PP^2$. Assuming $f$ is essential, it seems that the Euler number of the complement $U\defs\PP^2\stm Z$ coincides with $3\mi|{\rm Sing}\,Z|\mi\de$, where $\de\eq1$ if there is a line $l\sst Z$ with $|l\cap{\rm Sing}\,Z|\eq1$, and 0 otherwise (since the Euler number of $\C^*$ is 0). The number seems to be 1 if the coefficient of $z^d$ in $f$ is nonzero up to permutation of variables, and 0 otherwise.
\end{rem}

\begin{rem} \label{R8.4}
Let ${\bf x}\eq(x_1,\dots,x_n)$ be a coordinate system of $\C^n$ with $\dd_{\bf x}\eq(\dd_{x_1},\dots,\dd_{x_n})$ the dual vector fields to $\ddd\one {\bf x}\eq(\ddd\one x_1,\dots,\ddd\one x_n)$. A vector field $\xi$ of degree 0 on $\C^n$ is expressed as
$$\xi={\bf x}\one C\one ^t\!\dd_{\bf x},$$
for a matrix $C\eq(c_{i,j})$ ($c_{i,j}\ins\C$). Let ${\bf y}\eq(y_1,\dots,y_n)$ be another coordinate system with $A$ a matrix such that
$${\bf x}={\bf y}A,\q A\one ^t\!\dd_{\bf x}={}^t\!\dd_{\bf y},$$
using $\langle^t\!\dd_{\bf x},\ddd\one {\bf x}\rangle\eq{\rm Id}$, where ${\bf x}$ can be identified with $\ddd\one {\bf x}$. (One can use also the equality $\dd_{y_i}=\sum_j(\dd x_j/\dd y_i)\dd_{x_j}$, where $x_j\eq\msum_i\,y_ia_{i,j}$ and $\dd x_j/\dd y_i\eq a_{i,j}$.) We then get that
$$\xi={\bf x}\one C\one ^t\!\dd_{\bf x}={\bf y}AC\!A^{-1}\one ^t\!\dd_{\bf y}.$$
This implies for instance that the set of eigenvalues of $C$ is an invariant of a vector field $\xi$ of degree 0. Using this, we can show the following.
\end{rem}

\begin{prop} \label{P8.1}
Assume $\nu_{d+3}\nes0$ with $n\eq3$. Then a homogeneous polynomial $f$ is either extremely degenerated or a $\C$-linear combination of
\begin{equation} \label{8.8}
x^{d-2m}(y^2{-}2xz)^m\q(0\less m\less\tfrac{d}{2}),
\end{equation}
for some coordinate system $(x,y,z)$ of $\C^3$.
\end{prop}

\begin{proof}
The hypothesis means that the polynomial $f\eq\msum_{i,j,k}\,a_{i,j,k}x^iy^jz^k$ is annihilated by a non-zero vector field $\xi$ of degree 0. In the case the associated matrix $C$ is semi-simple, we may assume $\xi\eq ex\dd_x\pl e'y\dd_y\pl e''z\dd_z$, hence $f$ is extremely degenerated.
\sk
Assume $C$ has two Jordan blocks with eigenvalues $e,e'$, that is,
$$\xi=ex\dd_x\pl (x{+}ey)\dd_y\pl e'z\dd_z.$$
Let $F^x_p$ ($p\ins\Z$) be the increasing filtration on $\C[x,y,z]$ spanned over $\C$ by monomials $x^iy^jz^k$ with $i\less p$, and similarly for $F^y,F^z$. Taking $\Gr^{F^y}$, we verify that the coefficients $a_{i,j,k}$ must vanish unless $\ell(i,j,k)\defs ei{+}ej{+}e'k\eq0$, that is,
$${\rm supp}\,f\defs\{(i,j,k)\mid a_{i,j,k}\nes0\}\subset\ell^{-1}(0).$$
So $f$ is extremely degenerated except the case $e\eq e'\eq0$. In the last case (where $\xi\eq x\dd_y$), without taking $\Gr^{F^y}$, we see that ${\rm supp}\,f\sst\{j\eq0\}$, hence $f$ is extremely degenerated also in this case.
\sk
Assume finally $C$ has only one Jordan block with eigenvalue $e$, that is,
$$\xi\eq ex\dd_x\pl(x{+}ey)\dd_y\pl(y{+}ez)\dd_z.$$
Taking $\Gr^{F^y}\Gr^{F^z}$, we see that $e\eq0$, that is, $\xi\eq x\dd_y\pl y\dd_z$ (since $f\eq0$ if $e\nes0$). We then verify that the vector space consisting of homogeneous polynomials of degree $d$ annihilated by $\xi$ coincides with the vector subspace spanned over $\C$ by the polynomials in \eqref{8.8},
counting their dimensions. (Here we use the projection $\Z^3\ni(i,j,k)\mapsto(i,j)\in\Z^2$ and the ``restriction" by the inclusion of indices $\Z\ni i\mapsto(i,0)\in\Z^2$.) This finishes the proof of Proposition~\ref{P8.1}.
\end{proof}

\begin{prop} \label{P8.2}
Set
\begin{equation} \label{8.9}
g_m(x,y)\defs(x\pl y^2)^{d-2m}x^m\q(0\less m\less\tfrac{d}{2}),
\end{equation}
which are obtained by substituting $x\pl y^2$ and $\frac{1}{2}$ respectively into $x$ and $z$ in \eqref{8.8}. Assume $d\sgt4$. Then any $\C$-linear combination of these polynomials
$$g\defs\msum_{m=0}^{[d/2]}\,c_mg_m(x,y)\q(c_m\ins\C)$$
cannot be a weighted homogeneous reduced polynomial for any analytic local coordinates. 
\end{prop}

\begin{proof}
Consider first the case $d\eq2e$ with $e\ins\N$. Assume $g$ is a weighted homogeneous reduced polynomial for some analytic local coordinates. Since $g$ is reduced, the leading coefficient $c_0$ cannot vanish, and we may assume $c_0\eq1$. Set
$$P(t)\defs\msum_{m=0}^e\,c_mt^{e-m}\eq\mprod_{j=1}^e\,(t\mi\ga_j),$$
with $\ga_j\ins\C$. We then get
$$g\eq\mprod_{j=1}^e\,\bl((x{+}y^2)^2\mi\ga_jx\br),$$
(considering $x^{-e}g$).
Since $g$ is reduced, the complex numbers $\ga_j$ must be mutually distinct, hence the polynomial
$$\msum_{m=0}^e\,c_my^{4(e-m)}x^m$$
is reduced (substituting 0 into $y$ in $g$ and using the symmetry of binomial coefficients). We then get a contradiction applying Lemma~\ref{L8.1} below.
\sk
Assume now $d\eq2e{+}1$. By an argument similar to the above case, one can verify that $g$ is Newton nondegenerate, and the 0-dimensional faces of its Newton polygon consist of $(e{+}1,0)$, $(e,2)$, $(0,4e{+}2)$ (using the reducedness of $g$). This finishes the proof of Proposition~\ref{P8.2}.
\end{proof}

\begin{lem} \label{L8.1}
Let $g\eq\msum_{u\ges e}\,g_u$ be a semi-weighted-homogeneous polynomial in variables $x,y$ with weighted degree $\deg_{\bf w}g\gess e\ins\N$, where ${\rm wt}\,x\eq w_x\defs1$, ${\rm wt}\,y\eq w_y\defs\tfrac{1}{4}$, the $\bf w$-degree $u$ part $g_u$ vanishes for any $u\,{\notin}\,\tfrac{1}{2}\Z$, and $g_e$ is reduced. Assume $e\sgt2$, and $yg_{e+\frac{1}{2}}\eq\ga\one x\dd_yg_e$ for some $\ga\ins\C^*$. Then $g$ cannot be weighted homogeneous for any analytic local coordinates.
\end{lem}

\begin{proof}
Assume $\eta(g)\eq g$ for a local holomorphic vector field $\eta$. We have the decomposition
$$\eta\eq\msum_{u>-1}\,\eta_{{\bf w},u}$$
with $L_{\xi}\one\eta_{\bf w,u}\eq u\one\eta_{{\bf w},u}$, that is, $\eta_{{\bf w},u}$ has ${\bf w}$-degree $u\ins\Q$, where $L_{\xi}$ denotes the Lie derivation for the Euler field $\xi\defs w_x\one x\dd_x\pl w_y\one y\dd_y$.
\sk
We can verify that $\eta_{{\bf w},c}\eq0$ for any $c\slt0$ by increasing induction on $c\slt 0$. Indeed, if we have $\eta_{{\bf w},u}\eq0$ for any $u\slt c$, then we see that
$$\eta_{{\bf w},c}(g_e)\eq0$$
using the $\bf w$-weight decomposition of $\eta(g)\eq g$, which is compatible with the action of $\eta$ on $g$. We have the acyclicity (except for the top degree) of the Koszul complex for $\dd_xg_e,\dd_yg_e$, since $g_e$ has an isolated singularity at 0. It follows that
$$\eta_{{\bf w},c}=h\bl((\dd_yg_e)\dd_x\mi(\dd_xg_e)\dd_y\br),$$
for some $h\ins\C[x,y]$, whose $\bf w$-degree is equal to $-(e\mi1\mi\tfrac{1}{4})\slt0$ (since $e\sgt2$), hence
$$\eta_{{\bf w},c}\eq0\q\h{for any}\,\,\,c\slt0.$$
\sk
We then get
$$\eta_{{\bf w},0}(g_e)\eq g_e.$$
This is deduced from the equality $\eta(g)\eq g$ using the compatibility of the action with the $\bf w$-degree. Setting $\zeta\defs\eta_{{\bf w},0}\mi\tfrac{1}{e}\one\xi$, we then see that $\zeta(g_e)\eq0$. Using the acyclicity of the Koszul complex for $\dd_xg_e,\dd_yg_e$ explained above, this implies that
$$\zeta\eq0,\q\h{that is,}\q\eta_{{\bf w},0}\eq\tfrac{1}{e}\one\xi.$$
\sk
We now see that $\eta_{{\bf w},\frac{1}{2}}$ is a $\C$-linear combination of
$$xy^2\dd_x,\q y^6\dd_x,\q y^3\dd_y,$$
(where $xy^{-1}\dd_y$ cannot be included). From the equality $\eta(g)\eq g$, we deduce that
$$\eta_{{\bf w},0}\bl(g_{e+\frac{1}{2}}\br)+\eta_{{\bf w},\frac{1}{2}}(g_e)=g_{e+\frac{1}{2}},$$
with $\eta_{{\bf w},0}\bl(g_{e+\frac{1}{2}}\br)\eq(1\pl\tfrac{1}{2e})g_{e+\frac{1}{2}}$. We then get a nontrivial $\C$-linear relation among
$$xy^2\dd_x(g_e),\q y^6\dd_x(g_e),\q y^3\dd_y(g_e),\q g_{e+\frac{1}{2}}.$$
By hypothesis this implies a one among
$$xy^3\dd_x(g_e),\q y^7\dd_x(g_e),\q y^4\dd_y(g_e),\q x\dd_y(g_e),$$
which is equivalent to the inequality
$$\dim_{\C}(\dd g_e)^{e+\frac{3}{4}}<4.$$
Here $(\dd g_e)^{e+\frac{3}{4}}$ denotes the ${\bf w}$-degree $e{+}\frac{3}{4}$ part of the Jacobian ideal $(\dd g_e)$ for the weighted homogeneous polynomial $g_e$. This dimension however stays invariant under deformations of weighted homogeneous polynomials having isolated singularities with weights fixed (using \eqref{7.3}, see also \cite{JKSY0}), hence it can be calculated in the case the weighted homogeneous polynomial is written as $x^e\pl y^{4e}$, where the result is 4, since we assume $e\sgt2$. We thus get a contradiction. This finishes the proof of Lemma~\ref{L8.1}.
\end{proof}

\begin{rem} \label{R8.5}
In the case $d\eq3,4$ in Proposition~\ref{P8.2}, the polynomial $g$ has an isolated singularity of type $A_3$ and $A_7$ respectively assuming that $g$ is reduced in the case $d\eq4$. Their spectral numbers are $\frac{k}{4}$ with $k\ins\{3,4,5\}$ and $\frac{k}{8}$ with $k\ins\{5,\dots,11\}$ respectively. We can apply Theorem~\ref{T1} to these cases. We have the vanishing of $\chi(U)$ and $\de_{k/d}$ for $k\eq\tfrac{d}{2},\tfrac{3d}{2}$, hence \eqref{12} holds with $k''\eq{-}\infty$. (For $d\eq4$, the monodromy $(-1)$-eigenspace of the first Milnor fiber cohomology $H^1(F_{\!f},\C)_{-1}$ vanishes using \cite[Theorem 5.3]{kosz} or \cite{cons}.) Here we need \cite[Theorem 2.2]{bcm} to show the inclusion $d(\R_f\stm\R_Z)\sst[3,+\infty)$.
\end{rem}

\begin{rem} \label{R8.6}
It does not seem necessarily easy to determine $\R_f$ for extremely degenerated polynomials, for instance if $f\eq h$ or $hz$ with $h\defs x^4y^6\pl z^{10}$, where $\chi(U)\eq1$, $\de_{j/d}\eq\mu'_{j/d}\eq1$ for any $j\ins\N\cap(d,2d)$, and $\tfrac{j}{d}\ins(\R_Z\mi1)\stm\R_f$ for $j\eq3,4,5$.
(It seems rather difficult to find a good example with $\chi(U)\eq0$.)
\end{rem}

\part{Fundamental exact sequences of vanishing cycles} \label{Pa2}
\nin
In this part we study the fundamental exact sequences about the vanishing cycles, and prove Theorem~\ref{T10.1} which is a key to the proof of Theorem~\ref{T2}.

\section{Calculation of the Milnor fiber cohomology} \label{S9}
In the notation and assumption of the \nameref{intr}, set
$$n':=\dim Y=n{-}1,\q U:=Y\setminus Z,$$
with $j:U\into Y$ the canonical inclusion. Set
\begin{equation} \label{9.1}
\Lambda_d:=\bl\{\la\in\C^*\mid\la^d=1\br\}.
\end{equation}
Let $H^{\ssb}(\Ff,\C)_{\la}$, $H_c^{\ssb}(\Ff,\C)_{\la}$ be the $\la$-eigenspaces of the monodromy $T$. It is well known that $T^d=id$, and hence $H^{\ssb}(\Ff,\C)_{\la}=0$ for $\la\notin\Lambda_d$ (by using the geometric monodromy).
\sk
For $\la\in\Lambda_d$, we have the rank 1 local system $L_{\la}$ on $U$ such that
\begin{equation} \label{9.2}
H^{\ssb}(U,L_{\la})=H^{\ssb}(\Ff,\C)_{\la},\q H^{\ssb}_c(U,L_{\la})=H^{\ssb}_c(\Ff,\C)_{\la},
\end{equation}
see for instance \cite[3.1.1]{kosz}. More precisely, let $\pi:\Xt\to X$ be the blow-up at the origin, and set $\f:=\pi^*f$. The exceptional divisor $E$ is identified with $Y$. Define
\begin{equation} \label{9.3}
L_{\la}:=\psi_{\f,\la}\C_{\Xt}|_U.
\end{equation}
This is justified by \cite[Theorem 4.2]{BuSa1}, see also \cite[Section 1.6]{BuSa2}.
\sk
The second isomorphism of \eqref{9.2} follows from the first isomorphism by using the self-duality of the nearby cycle sheaves $\psi_{\f}\C_{\Xt}$, which implies the isomorphisms
\begin{equation} \label{9.4}
L_{\la}^{\vee}=L_{\la^{-1}}\,\,\,(\la\in\Lambda_d),\q\q L_1=\C_U,
\end{equation}
see also \cite[Sections 1.3--4]{BuSa2}. The shifted local system $\mopl_{\la\in\Lambda_d}L_{\la}[n']$ naturally underlie a mixed Hodge module on $U$ by \eqref{9.3}, and the isomorphisms in \eqref{9.2} are compatible with mixed Hodge structures.
\sk
Note that the rank 1 local systems $L_{\la}$ are uniquely determined by the local monodromies at smooth points of $Z$, which are given by the multiplication by $\la^{-1}$. Indeed, in the notation of \eqref{4.1}, $L_{\la}$ is uniquely determined by its restriction over $U\cap Y'$ (using the direct image by $U\cap Y'\into U$), and we have the isomorphism
\begin{equation} \label{9.5}
L_{\la}|_{U\cap Y'}=h^*E_{\la^{-1}},
\end{equation}
where $E_{\la^{-1}}$ is a rank 1 local system on $\C^*$ with monodromy $\la^{-1}$. (We can verify that $L_{\la}|_{U\cap Y'}\otimes h^*E_{\la}$ is extendable as a rank 1 local system over $Y'\setminus\Si$, and then over $Y'=\C^{n'}$ since $\dim Y'\ges 2$.)
\sk
Note that \eqref{9.5} also implies \eqref{9.4} as well as the isomorphisms
\begin{equation} \label{9.6}
\psi_{h,\eta}L_{\la}|_{U\cap Y'}=\psi_{h,\la\eta}\C_{U\cap Y'}\q\q(\forall\,\eta\in\C^*).
\end{equation}
Indeed, we have more generally
\begin{equation} \label{9.7}
\psi_{h,\la}\F=\psi_{h,1}(\F\otimes h^*E_{\la^{-1}})\q\h{for any}\,\,\,\F\in D^b_c(\C_{U\cap Y'}).
\end{equation}
\sk
It is well known (see \cite{Mi} and also \cite[2.1.2]{cons}, etc.) that if $\dim{\rm Sing}\,f^{-1}(0)=1$, then
\begin{equation} \label{9.8}
\Ht^j(U,L_{\la})=\Ht^{\ssb}(\Ff,\C)_{\la}=0\q\h{unless}\,\,\,j\in\{n'{-}1,n'\},
\end{equation}
with $\Ht^j(U,L_{\la}):=H^j(U,L_{\la})$ for $\la\ne 1$.
(Recall that $L_1=\C_U$, see \eqref{9.4}.)
\sk
For $z\in\Si={\rm Sing}\,Z$ (see \eqref{3.1}), there are distinguished triangles
\begin{equation} \label{9.9}
(\Rb j_*L_{\la})_z\to\psi_{h,1}(L_{\la}|_{U\cap Y'})_z\buildrel{N}\over\to\psi_{h,1}(L_{\la}|_{U\cap Y'})_z(-1)\buildrel{+1}\over\to.
\end{equation}
The associated long exact sequence is strictly compatible with the Hodge filtration $F$ coming from the theory of mixed Hodge modules \cite{mhm} (by using for instance \cite[Proposition 1.3]{ext}).
\sk
By \eqref{9.4} for $\eta=1$ together with condition~(WH) in Theorem~\ref{T2}, we have the vanishing of $N$ in the long exact sequences associated with the distinguished triangles in \eqref{9.9}. Using this, we can prove the canonical isomorphisms for $z\in\Si$
\begin{equation} \label{9.10}
(R^ij_*L_k)_z\cong\begin{cases}H^{n'-1}(\Fhz)_{\la}&\h{if}\,\,\,i=n'{-}1,\\ H^{n'-1}(\Fhz)_{\la}(-1)&\h{if}\,\,\,i=n',\\ \,0&\h{if}\,\,\,i\in[2,n'{-}2]\,\,\,\h{or}\,\,\,i>n',
\end{cases}
\end{equation}
where $(p)$ for $p\in\Z$ denotes the Tate twist in general (see \cite{th2}), $\Fhz$ is the Milnor fiber of $h_z$, and the coefficient field $\C$ of the cohomology is omitted to simplify the notation. These isomorphisms are compatible with the Hodge filtration $F$ and the weight filtration $W$ of mixed Hodge structures.
\sk
Actually we can also construct the isomorphisms in \eqref{9.10} by using another method. Indeed, let $\Q_{h,X}:=a_X^*\Q\in D^b{\rm MHM}(X)$ with $a_X:X\to pt$ the structure morphism. Let $i_0:\{0\}\into X:=\C^n$, $j_0:X\setminus\{0\}\into X$ be natural inclusions. Using \cite[Theorem 4.2]{BuSa1} together with the distinguished triangle
$$(i_0)_*i_0^!\varphi_f\Q_{h,X}\to \varphi_f\Q_{h,X}\to(j_0)_*j_0^*\varphi_f\Q_{h,X}\buildrel{+1}\over\to,$$
we can prove Theorem~\ref{T10.1} below by showing the self-dual exact sequences of mixed Hodge structures (see also \cite[2.1]{cons}):
\begin{equation} \label{9.11}
\aligned0&\to H^{n'-1}(\Ff)_{\ne 1}\to\mopl_{z\in\Si}\,H^{n'-1}(\Fhz)^{T_l}_{\ne 1}\buildrel{\rho^{\vee}}\over\to H_c^{n'}(\Ff)_{\ne 1}\\&\to H^{n'}(\Ff)_{\ne 1}\buildrel{\rho}\over\to\mopl_{z\in\Si}\,H^{n'-1}(\Fhz)^{T_l}_{\ne 1}(-1)\to H_c^{n'+1}(\Ff)_{\ne 1}\to 0,\endaligned
\end{equation}
\begin{equation} \label{9.12}
\aligned0&\to H^{n'-1}(\Ff)_1\to\mopl_{z\in\Si}\,H^{n'-1}(\Fhz)^{T_l}_1\buildrel{\rho^{\vee}}\over\to H_c^{n'}(\Ff)_1(-1)\\&\to H^{n'}(\Ff)_1\buildrel{\rho}\over\to\mopl_{z\in\Si}\,H^{n'-1}(\Fhz)^{T_l}_1(-1)\to H_c^{n'+1}(\Ff)_1(-1)\to 0,\endaligned
\end{equation}
where $E^{\,T_l}$ denotes the invariant part by the action of the local system monodromy $T_l$ on $E:=H^{n'-1}(\Fhz)$.
These are called the {\it fundamental exact sequences of vanishing cycles}.
Since this proof is rather long, it is omitted in this version.
In the proof of Theorem~\ref{T10.1} in Sections~\ref{S11}--\ref{S12} below, we will construct long exact sequence which should correspond to the above ones (except for the case $\la=1$, $n\eq3$). For the moment it is {\it still unclear} whether these two canonical isomorphisms {\it really coincide} (even though this is very much expected to hold).

\section{Key theorem to the proof of Theorem~\ref{T2}} \label{S10}
We show the following theorem whose proof will be given in Sections~\ref{S11}--\ref{S14}.

\begin{thm} \label{T10.1}
With the notation and assumption of Section~{\rm\ref{S9}} together with {\rm (WH)} in Theorem~$\ref{T2}$, there is a canonical isomorphism compatible with the Hodge filtration $F:$
\begin{equation} \label{10.1}
{\rm Coker}\Bigl(H^{n'}(\Ff)_{\la}\buildrel{\rho}\over\to\bigoplus_{z\in\Si}\,H^{n'-1}(\Fhz)_{\la}(-1)\Bigr)=\begin{cases}H^{n'-1}(\Ff)_{\la^{-1}}^{\vee}(-n')\raise-3mm\h{\,}&\h{if}\,\,\,\la\ne 1,\\H^{n'-1}(\Ff)_1^{\vee}(-n'{-}1)&\h{if}\,\,\,\la=1,\end{cases}
\end{equation}
where $\Ff$ is the Milnor fiber of $f$ as in the \nameref{intr}, $\rho$ is a canonical morphism, and the coefficient field $\C$ of the cohomology is omitted to simplify the notation. Moreover $\rho$ can be identified via the isomorphisms as in {\rm\eqref{9.2}} and {\rm\eqref{9.10}} with the following canonical morphism\,:
\begin{equation} \label{10.2}
H^{n'}(U,L_{\la})\to\mopl_{z\in\Si}\,(R^{\one n'}j_*L_k)_z.
\end{equation}
\end{thm}

\begin{rem} \label{R10.1}
By the first injective morphisms in \eqref{9.11}--\eqref{9.12}, $H^{n'-1}(\Ff)_{\ne 1}$ and $H^{n'-1}(\Ff)_1$ are identified respectively with subspaces of
$$\mopl_{z\in\Si}\,H^{n'-1}(\Fhz)_{\ne 1},\q\mopl_{z\in\Si}\,H^{n'-1}(\Fhz)_1,$$
which are {\it pure $\Q$-Hodge structures of weight $n'{-}1$ and $n'$} under the assumption~(WH) in Theorem~\ref{T2}, see \cite{St1}. In particular, polarizations of Hodge structures induce self-duality isomorphisms of $\Q$-Hodge structures
\begin{equation} \label{10.3}
H^{n'-1}(\Ff)_{\ne 1}\cong H^{n'-1}(\Ff)_{\ne 1}^{\vee}(1-n'),\q H^{n'-1}(\Ff)_1\cong H^{n'-1}(\Ff)_1^{\vee}(-n').
\end{equation}
These are closely related to the difference in the Tate twist on the right-hand side of \eqref{10.1} for $\la\ne 1$ and $\la=1$.
\end{rem}

We give a proof of Theorem~\ref{T10.1} (which is closely related to \cite{Di1}, \cite{Kl}) by dividing it into three cases.

\section{Proof of Theorem~\ref{T10.1}, Part I} \label{S11}
We give the proof for the case $\la\ne 1$ in this section. Set
$$Y^{\circ}:=Y\setminus\Si\,\,\,\,\h{with}\,\,\,j':Y^{\circ}\into Y\,\,\,\h{the inclusion.}$$
Let $L'_{\la}$ be the zero extension of $L_{\la}$ over $Y^{\circ}$. There is a distinguished triangle
\begin{equation} \label{11.1}
j'_!L'_{\la}\to\Rb j'_*L'_{\la}\to\mopl_{z\in\Si}\,(\Rb j_*L_{\la})_z\buildrel{+1}\over\to,
\end{equation}
together with isomorphisms
\begin{equation} \label{11.2}
j'_!L'_{\la}=j_!L_{\la},\q\Rb j'_*L'_{\la}=\Rb j_*L_{\la}.
\end{equation}
These induce exact sequences
\begin{equation} \label{11.3}
\aligned0&\to H^{n'-1}(U,L_{\la})\to\mopl_{z\in\Si}\,(R^{\one n'-1}j_*L_{\la})_z\to H_c^{n'}(U,L_{\la})\\&\to H^{n'}(U,L_{\la})\to\mopl_{z\in\Si}\,(R^{\one n'}j_*L_{\la})_z\to H_c^{n'+1}(U,L_{\la})\to 0,\endaligned
\end{equation}
which are essentially self-dual by replacing $\la$ with $\la^{-1}$.
\sk
By \eqref{9.2}, \eqref{9.10}, these give essentially self-dual exact sequences
\begin{equation} \label{11.4}
\aligned0&\to H^{n'-1}(\Ff)_{\la}\to\mopl_{z\in\Si}\,H^{n'-1}(\Fhz)_{\la}\buildrel{\rho^{\vee}}\over\to H_c^{n'}(\Ff)_{\la}\\&\to H^{n'}(\Ff)_{\la}\buildrel{\rho}\over\to\mopl_{z\in\Si}\,H^{n'-1}(\Fhz)_{\la}(-1)\to H_c^{n'+1}(\Ff)_{\la}\to 0.\endaligned
\end{equation}
It is expected that their direct sum over $\la\in\Lambda_d\setminus\{1\}$ would be identified with \eqref{9.11}. Applying Poincar\'e duality to the last term of \eqref{11.4}, we get Theorem~\ref{T10.1} in the case $\la\ne 1$.

\section{Proof of Theorem~\ref{T10.1}, Part II} \label{S12}
We give the proof for the case $\la=1$, $n'\ges 3$ in this section. 
We assume for the moment {\it only} $\la=1$. Set
$$Z^{\circ}:=Z\setminus\Si\,\,\,\,\h{with}\,\,\,j_Z:Z^{\circ}\into Z\,\,\,\h{the inclusion.}$$
We have the distinguished triangle
\begin{equation} \label{12.1}
\Q_{Y^{\circ}}\to\Rb j''_*\Q_U\to\Q_{Z^{\circ}}(-1)[-1]\buildrel{+1}\over\to,
\end{equation}
where $j'':U\into Y^{\circ}$ is the inclusion, and $Y^{\circ}:=Y\setminus\Si$ with $j':Y^{\circ}\into Y$ as in Section~\ref{S11}. Applying $\Rb j'_*$ to \eqref{12.1}, we get the distinguished triangle
\begin{equation} \label{12.2}
\Rb j'_*\Q_{Y^{\circ}}\to\Rb j_*\Q_U\to\Rb(j_Z)_*\Q_{Z^{\circ}}(-1)[-1]\buildrel{+1}\over\to.
\end{equation}
These naturally underlie distinguished triangles of complexes of mixed Hodge modules.
\sk
Setting $m:=|\Si|$, we have
\begin{equation} \label{12.3}
H^i(Y^{\circ})=\begin{cases}\,\Q(-j)&\h{if}\,\,\,i=2j\,\,\,\h{with}\,\,\,j\in[0,n'{-}1],\\ \buildrel{m-1}\over\mopl\Q(-n')&\h{if}\,\,\,i=2n'{-}1,\\ \,\,\,0&\h{otherwise}.\end{cases}
\end{equation}
\begin{equation} \label{12.4}
\Ht^i(U)=0\q\h{unless}\,\,\,\,i\in\{n'{-}1,\,n'\},\,\,\h{see (2.1.8)}.
\end{equation}
\sk
{\it Assume $n'\ges 3$ from now on.} The distinguished triangle \eqref{12.2} implies the isomorphisms for $z\in\Si:$
\begin{equation} \label{12.5}
(R^ij_*\Q_U)_z=(R^{\one i-1}(j_Z)_*\Q_{Z^{\circ}})_z(-1)\q\q(i\in[1,2n'{-}3]),
\end{equation}
as well as the short exact sequences of mixed Hodge structures
\begin{equation} \label{12.6}
0\to H^{i+1}(U)\to H^i(Z^{\circ})(-1)\to H^{i+2}(Y^{\circ})\to 0\q(i\in\N),
\end{equation}
where the coefficient field $\Q$ of the cohomology is omitted to simplify the notation. Indeed, the restriction morphisms $H^i(Y^{\circ})\to H^i(U)$ vanish for any $i>0$ (where we use \eqref{12.4} together with $2n'{-}1>n'$ in the case $i=2n'{-}1$ in \eqref{12.3}).
\sk
Define the {\it primitive part} by
$$H^i_{\rm prim}(Z^{\circ}):={\rm Ker}\bl(H^i(Z^{\circ})\to H^{i+2}(Y^{\circ})(1)\br)\q\h{for}\,\,\,\,i\in\{n'{-}2,n'{-}1\},$$
where the morphism on the right-hand side is the Gysin morphism appearing in \eqref{12.6}. (This is compatible with the standard definition of the primitive part.) Note that the same definition does not work very well for $i=n'$ if $n'=3$, see \eqref{12.3}. We are only interested in the case $i\in\{n'{-}2,n'{-}1\}$ for the exact sequence \eqref{12.10} explained below.
\sk
By \eqref{12.6} we get the isomorphisms
\begin{equation} \label{12.7}
H^i_{\rm prim}(Z^{\circ})=H^{i+1}(U)(1)\q\q\bl(i\in\{n'{-}2,n'{-}1\}\br).
\end{equation}
\sk
Define the {\it primitive part} for the cohomology with compact supports by
$$H^i_{c,{\rm prim}}(Z^{\circ}):=H^{2n'{-}2-i}_{\rm prim}(Z^{\circ})^{\vee}(1-n')\q\q\bl(i\in\{n'{-}1,n'\}\br).$$
By \eqref{12.7} this vanishes for $i=n'{-}2$, and we have
\begin{equation} \label{12.8}
H^i_{c,{\rm prim}}(Z^{\circ})=H^{2n'-1-i}(U)^{\vee}(-n')=H_c^{i+1}(U)\q\q\bl(i\in\{n'{-}1,n'\}\br).
\end{equation}

On the other hand, there is a distinguished triangle
\begin{equation} \label{12.9}
(j_Z)_!\Q_{Z^{\circ}}\to\Rb(j_Z)_*\Q_{Z^{\circ}}\to(\Rb(j_Z)_*\Q_{Z^{\circ}})|_{\Si}\buildrel{+1}\over\to,
\end{equation}
underlying naturally a distinguished triangle of complexes of mixed Hodge modules.
By Section~\ref{S13} below, this induces a self-dual exact sequence of mixed Hodge structures
\begin{equation} \label{12.10}
\aligned 0&\to H^{n'-2}_{\rm prim}(Z^{\circ})\to\mopl_{z\in\Si}\,(R^{\one n'-2}(j_Z)_*\Q_{Z^{\circ}})_z\to H_{c,{\rm prim}}^{n'-1}(Z^{\circ})\\&\to H_{\rm prim}^{n'-1}(Z^{\circ})\to\mopl_{z\in\Si}\,(R^{\one n'-1}(j_Z)_*\Q_{Z^{\circ}})_z\to H_{c,{\rm prim}}^{n'}(Z^{\circ})\to 0,\endaligned
\end{equation}
assuming $n'\ges 3$. Combined with \eqref{9.10}, \eqref{12.5}, and \eqref{12.7}--\eqref{12.8}, this gives a self-dual exact sequence of mixed Hodge structures
\begin{equation} \label{12.11}
\aligned 0&\to H^{n'-1}(U)\to\mopl_{z\in\Si}\,H^{n'-1}(\Fhz)_1\to H_c^{n'}(U)(-1)\\&\to H^{n'}(U)\to\mopl_{z\in\Si}\,H^{n'-1}(\Fhz)_1(-1)\to H_c^{n'+1}(U)(-1)\to 0,\endaligned
\end{equation}
(which is expected to be identified with \eqref{9.12} in this case). Theorem~\ref{T10.1} then follows in the case $\la=1$ and $n'\ges 3$ by using Poincar\'e duality for $U$ together with \eqref{9.2}.

\section{Proof of the exact sequence (\ref{12.10})} \label{S13}
Define the {\it non-primitive part} by
$$\aligned H^i_{c,{np}}(Z^{\circ})&:={\rm Im}\bl(H_c^i(Y^{\circ})\to H_c^i(Z^{\circ})\br)\q\q\bl(i\in\{n'{-}1,n'\}\br),\\ H^i_{np}(Z^{\circ})&:={\rm Im}\bl(H^i(Y^{\circ})\to H^i(Z^{\circ})\br)\q\q\bl(i\in\{n'{-}1,n'{-}2\}\br),\endaligned$$
\sk
For $i=n'{-}1$, we have a canonical isomorphism
\begin{equation} \label{13.1}
H^{n'-1}_{c,{np}}(Z^{\circ})\simto H^{n'-1}_{np}(Z^{\circ}),
\end{equation}
by using the canonical isomorphism $H^{n'-1}_c(Y^{\circ})\simto H^{n'-1}(Y^{\circ})$.
\sk
For $i=n'$, a similar argument implies the injectivity of the canonical isomorphism
\begin{equation} \label{13.2}
H^{n'}_{c,{np}}(Z^{\circ})\into H^{n'}(Z^{\circ}).
\end{equation}
\sk
There are moreover canonical isomorphisms
\begin{equation} \label{13.3}
\begin{array}{rll}H^i_{c,{\rm prim}}(Z^{\circ})=&\!\!\!H_c^i(Z^{\circ})/H^i_{c,{np}}(Z^{\circ})&\q\bl(i\in\{n'{-}1,n'\}\br),\raise-3mm\h{}\\ H^i_{\rm prim}(Z^{\circ})=&\!\!\!H^i(Z^{\circ})/H^i_{np}(Z^{\circ})&\q\bl(i\in\{n'{-}1,n'{-}2\}\br).\end{array}
\end{equation}
Indeed, the first isomorphism easily follows from the definition (since the dual of the kernel of a morphism is the cokernel of the dual morphism in general), and the second follows from the bijectivity of the composition
$$H^i(Y^{\circ})\to H^i(Z^{\circ})\to H^{i+2}(Y^{\circ})(1)\q\q\bl(i\in\{n'{-}1,n'{-}2\}\br).$$
\sk
By \eqref{13.1} and \eqref{13.3} for $i=n'{-}1$, we get the exactness of the middle part of \eqref{12.10} using the long exact sequence associated with \eqref{12.9}. By the self-duality of \eqref{12.10} together with the first isomorphism of \eqref{13.3} for $i=n'$ and the injectivity of \eqref{13.2}, it now remains to show the surjectivity of the morphism to $H_{c,{\rm prim}}^{n'}(Z^{\circ})$ in \eqref{12.10}.
\sk
Since $H_{c,{\rm prim}}^{n'}(Z^{\circ})$ has weights $\les n'$ (see \cite{th2}, \cite{mhm}), it is enough to show
\begin{equation} \label{13.4}
H^{n'}(Z^{\circ})/H^{n'}_{c,{np}}(Z^{\circ})\,\,\,\h{has weights}>n',
\end{equation}
by using the long exact sequence associated with \eqref{12.9} together with the injectivity of \eqref{13.2}. However, \eqref{13.4} easily follows from \eqref{12.3}--\eqref{12.4} and \eqref{12.6}. This finishes the proof of the exact sequence \eqref{12.10}.

\section{Proof of Theorem~\ref{T10.1}, Part III} \label{S14}
We give the proof for the case $\la=1$, $n'=2$.
Here $Z$ is a curve on $Y=\PP^2$. Using the distinguished triangle \eqref{12.2} together with the snake lemma, we can get the commutative diagram of exact sequences
\begin{equation} \label{14.1}
\begin{gathered}
\xymatrix{H^2(U)\,\ar@{^(->}[r]\ar[d]^{\rho}&H^1(Z^{\circ})(-1)\ar@{->>}[r]\ar[d]&H^3(Y^{\circ})\ar@{^(->}[d]\\ \mopl_z\,(R^2j_*\Q_U)_z\,\ar@{^(->}[r]\ar@{->>}[d]&\mopl_z\,(R^1(j_Z)_*\Q_{Z^{\circ}}(-1))_z\ar@{->>}[r]\ar@{->>}[d]&\mopl_z\,(R^3j'_*\Q_{Y^{\circ}})_z\ar@{->>}[d]\\ {\rm Coker}\,\rho\,\ar@{^(->}[r]&H_c^2(Z^{\circ})(-1)\ar@{->>}[r]&H_c^4(Y^{\circ})}
\end{gathered}
\end{equation}
where the direct sums are taken over $z\in\Si$, see \cite[Ch.~6, 3.14]{Di1} and Remark~\ref{R14.1} below. (As for the right vertical exact sequence, note that $\dim H^3(Y^{\circ})=m-1$, $\dim H_c^4(Y^{\circ})=1$.)
\sk
We have moreover the canonical isomorphism
\begin{equation} \label{14.2}
{\rm Coker}\,\rho=H^1(U)^{\vee}(-3),
\end{equation}
using the dual of the isomorphism
\begin{equation} \label{14.3}
H^1(U)={\rm Coker}\bl(H^0(Y^{\circ})\to H^0(Z^{\circ})\br)(-1).
\end{equation}
The latter follows from \eqref{12.6} for $i=0$ by using the bijectivity of the composition
$$H^0(Y^{\circ})\to H^0(Z^{\circ})\to H^2(Y^{\circ})(1).$$
So \eqref{10.1} also follows in the case $\la=1$, $n'=2$. This finishes the proof of Theorem~\ref{T10.1}.

\begin{rem} \label{R14.1}
In the case $n'=2$, let $r_z$ and $r_Z$ be respectively the number of local and global irreducible components of $(Z,z)$ and $Z$. Then
\begin{equation} \label{14.4}
\dim(R^2j_*\Q_U)_z=r_z-1,\q\dim(R^1(j_Z)_*\Q_{Z^{\circ}})_z=r_z,\q\dim(R^3j'_*\Q_{Y^{\circ}})_z=1.
\end{equation}
\vskip-7mm
\begin{equation} \label{14.5}
\dim H^1(U)=r_Z-1,\q\dim H^0(Z^{\circ})=r_Z,\q\dim H^0(Y^{\circ})=\dim H^2(Y^{\circ})=1.
\end{equation}
\end{rem}

\section{Relation with the Euler characteristics} \label{S15}
Set $U:=Y\setminus Z$. By \eqref{9.2} we have
\begin{equation} \label{15.1}
\msum_i\,(-1)^i\dim H^i(\Ff,\C)_{\la}=\chi(U)\q\bl(\la\in\Lambda_d\br).
\end{equation}
In the case the pole order spectral sequence {\it degenerates at} $E_2$ as in Theorem~\ref{T2}, we then get
\begin{equation} \label{15.2}
\msum_{j\in\Z}\,\bl(\mu_{k+jd}-\nu_{k+jd+d}\br)=\begin{cases}(-1)^{n'}\chi(U)&\h{if}\,\,\,k\in\Z\setminus d\,\Z,\\ (-1)^{n'}\bl(\chi(U)-1\br)&\h{if}\,\,\,k\in d\,\Z.\end{cases}
\end{equation}
Here $n'=n{-}1$. On the other hand, it is well known that $\chi(U)=n'+1-\chi(Z)$ (since $\chi(\PP^{n'})=n'+1$), and for a {\it smooth} hypersurface $V\subset\PP^{n'}$ of degree $d$, we have
\begin{equation} \label{15.3}
\chi(Z)=\chi(V)-(-1)^{n'-1}\mu_Z\q\h{with}\q\mu_Z:=\msum_{z\in\Si}\,\mu_z,
\end{equation}
by using a one-parameter family $\{f\,{+}\,sg\}_{s\in\C}$ with $V=\{g=0\}$, where $\mu_z$ is the Milnor number of $h_z$, see \cite{van}. (Here we may assume $V\cap{\rm Sing}\,Z=\emptyset$ since ${\rm Sing}\,Z$ is isolated.) So we get
\begin{equation} \label{15.4}
\chi(U)=\chi(\PP^{n'}\setminus V)-(-1)^{n'}\mu_Z.
\end{equation}
\sk
In the case $n'=2$, we have
\begin{equation} \label{15.5}
\chi(U)=(d-1)(d-2)+1-\mu_Z.
\end{equation}
since a well-known formula for the genus of smooth plane curves implies that
\begin{equation} \label{15.6}
\chi(\PP^2\setminus V)=(d-1)(d-2)+1.
\end{equation}

\part{Calculation of twisted de Rham complexes} \label{Pa3}
\nin
In this part we calculate the filtered twisted de Rham complexes of weighted homogeneous polynomials having isolated singularities.

\section{Twisted de Rham complexes} \label{S16}
Let $h$ be a weighted homogeneous polynomial for a local coordinate system $(y_1,\dots,y_{n'})$ of $(Y,0):=(\C^{n'},0)$ (as a complex manifold) with positive weights $w_1,\dots,w_{n'}$, see Remark~\ref{R4.1}. We assume moreover that $h^{-1}(0)$ has an {\it isolated singularity} at 0. We take and fix
$$\al\in\Q\cap(0,1].$$
Consider the {\it twisted de Rham complex}
\begin{equation} \label{16.1}
K^{\ssb}(h,\al):=\Om_{Y,0}^{\ssb}[h^{-1}]h^{-\al}\cong\bl(\Om_{Y,0}^{\ssb}[h^{-1}];\ddd-\al\tfrac{\ddd h}{h}\wedge\br),
\end{equation}
where $\Om_Y^{\ssb}$ is the analytic de Rham complex.
Define the {\it pole order filtration} $P$ by
\begin{equation} \label{16.2}
P^kK^i(h,\al):=\Om_{Y,0}^i\,h^{-\al+k-i}\q\q(i,k\in\Z),
\end{equation}
so that
\begin{equation} \label{16.3}
\Gr_P^kK^{\ssb}(h,\al)\cong(\Om_{Y,0}^{\ssb}/h\Om_{Y,0}^{\ssb};\ddd h\wedge).
\end{equation}
Set $P_k=P^{-k}$ as usual. Note that $K^{\ssb}(h,\al)$ is a {\it graded complex} with degrees defined by
\begin{equation} \label{16.4}
\deg y_i=\deg\ddd y_i=w_i,\q\deg h=\deg\ddd h=1,\q\deg h^{-\al}=-\al,
\end{equation}
and $P$ is compatible with the grading. Note that the $H^i\bl(K^{\ssb}(h,\al)\br)$ ($i\in\Z$) are graded $\C$-vector spaces, and have the induced filtration $P$ compatible with the grading.
\sk
Set $\la:=\exp(-2\pi i\al)$. We have the isomorphism
\begin{equation} \label{16.5}
H^{n'}\!\bl(K^{\ssb}(h,\al)\br)=H^{n'-1}(F_h,\C)_{\la}(-1).
\end{equation}
This follows from a generalization of \eqref{9.10} with $L_{\la}$ replaced by $h^*E_{\la^{-1}}$ as in \eqref{9.5}, where $h$ is a weighted homogeneous polynomial for some local coordinates of a complex manifold $Y'$. Indeed, $\OO_Y[h^{-1}]h^{-\al}$ is a regular holonomic $\D_Y$-module corresponding to $\Rb j_*h^*E_{\la^{-1}}[n']$ by the de Rham functor ${\rm DR}_Y$ so that we have a canonical isomorphism in the derived category
\begin{equation} \label{16.6}
K^{\ssb}(h,\al)=\Rb j_*h^*E_{\la^{-1}},
\end{equation}
where $E_{\la^{-1}}$ is as in the explanation after \eqref{9.5}, and $j$ is the inclusion of the complement of $h^{-1}(0)$.

\section{Brieskorn modules} \label{S17}
In the above notation and assumption, define the Brieskorn module $H''_h$ (see \cite{Br}) by
\begin{equation} \label{17.1}
H''_h:=\Om_{Y,0}^{n'}/\ddd h\,{\wedge}\,\ddd\Om_{Y,0}^{n'-2},
\end{equation}
where $n'=\dim Y$. This is a {\it graded module} with degrees defined by \eqref{16.4}.
\sk
The actions of $t$ and $\dd_tt$ on $H''_h$ are defined by
\begin{equation} \label{17.2}
t[\om]:=[h\om],\q\dd_tt[\om]:=\ddd\eta\q\h{with}\q\ddd h\,{\wedge}\,\eta=h\om,
\end{equation}
where $\om\in\Om_{Y,0}^{n'}$, $\eta\in\Om_{Y,0}^{n'}$. It is well known that $H''_h$ is a free $\C\{t\}$-module of rank $\mu$ with $\mu$ the Milnor number of $h$. For the definition of $\dd_tt$, we use a property of the weighted homogeneous polynomial that $h$ is contained in the Jacobian ideal $(\dd h)\subset\OO_{Y,0}$, see Remark~\ref{R4.1}. More precisely, using \eqref{4.8}, it is easy to show that
\begin{equation} \label{17.3}
\dd_tt[\om]=\beta[\om]\q\h{for}\,\,\,[\om]\in H''_{h,\beta},
\end{equation}
(see the proof of Proposition~\ref{P18.1} below), where $H''_{h,\beta}\subset H''_h$ is the degree $\beta$ part of graded module. (It is then easy to show the $t$-torsion-freeness of $H''_h$ by using \eqref{17.3} in the weighted homogeneous polynomial case, since $H''_{h,\beta}=0$ for $\beta\les 0$.)

\section{Relation with Brieskorn modules} \label{S18}
In this section we prove the following.

\begin{prop} \label{P18.1}
With the above notation and assumption $($in particular, $\al\in(0,1])$, we have
\begin{equation} \label{18.1}
\aligned H^{n'}\!\bl(K^{\ssb}(h,\al)\br)_{\beta}&=0\q(\forall\,\beta\ne 0),\q\h{that is,}\\ H^{n'}\!\bl(K^{\ssb}(h,\al)\br)&=H^{n'}\!\bl(K^{\ssb}(h,\al)\br)_0,\endaligned
\end{equation}
and there are canonical isomorphisms
\begin{equation} \label{18.2}
\iota_k:H''_{h,\al+k}\simto P_{k-n'}H^{n'}\!\bl(K^{\ssb}(h,\al)\br)\q\q(k\in\N),
\end{equation}
satisfying
\begin{equation} \label{18.3}
\iota_{k+i}\ssc t^i=\iota_k\q\h{in}\,\,\,H^{n'}\!\bl(K^{\ssb}(h,\al)\br)\q\q(k,i\in\N),
\end{equation}
where $t^i:H''_{h,\al+k}\into H''_{h,\al+k+i}$ is induced by the action of $t$ on $H''_h$ defined in {\rm\eqref{17.2}}.
\end{prop}

\begin{proof}
Set
$$\om_0:=\ddd y_1\,{\wedge}\,\cdots\,{\wedge}\,\ddd y_{n'},\q\eta_0:=\msum_i\,(-1)^{i-1}w_iy_i\,\ddd y_1\,{\wedge}\,\cdots\,{\wedge}\,\widehat{\ddd y_i}\,{\wedge}\,\cdots\,{\wedge}\,\ddd y_{n'}.$$
The latter is the inner derivation (or contraction) of $\om_0$ by the vector field $v$ in \eqref{4.8} so that
\begin{equation} \label{18.4}
\ddd h\wedge g\,\eta_0=hg\,\om_0\q\h{for}\,\,\,g\in\OO_{Y,0}.
\end{equation}
Set $\beta_0:=\msum_iw_i=\deg\om_0$. We can easily verify that
\begin{equation} \label{18.5}
\ddd(g\,\eta_0)=\beta g\,\om_0\q\h{for}\,\,\,g\in(\OO_{Y,0})_{\beta-\beta_0}\,\,\,\bl(\h{that is,}\,\,\, g\om_0\in(\Om_{Y,0}^{n'})_{\beta}\br),
\end{equation}
by using the equality $v(g)=\beta g$ for $g\in(\OO_{Y,0})_{\beta}$ together with $\msum_iw_i\dd_{y_i}y_i=v+\beta_0$.
\sk
The assertion \eqref{18.1} then follows from \eqref{18.4}--\eqref{18.5} by the definition of the differential of $K^{\ssb}(h,\al)$ in \eqref{16.1}. (Note that \eqref{17.3} also follows from \eqref{18.4}--\eqref{18.5}.)
\sk
Define the canonical morphisms $\iota_k$ in \eqref{18.2} by
\begin{equation} \label{18.6}
\iota_k[\om]=[\om h^{-\al-k}]\in H^{n'}\!\bl(K^{\ssb}(h,\al)\br)_0\q\h{for}\,\,\,\om\in(\Om_{Y,0}^{n'})_{\al+k}.
\end{equation}
Then \eqref{18.3} and the surjectivity of $\iota_k$ in \eqref{18.2} follow from the definition. For the well-definedness of \eqref{18.6}, we have to show that
\begin{equation} \label{18.7}
\iota_k[\om]=0\q\h{if}\q\om=\ddd h\,{\wedge}\,\ddd\xi\q\h{with}\q\xi\in(\Om_{Y,0}^{n'-2})_{\al+k-1}.
\end{equation}
Set $\eta:=\ddd\xi$. Then
\begin{equation} \label{18.8}
\ddd(\eta\,h^{-\al-k+1})=-(\al+k-1)\ddd h\wedge\eta=-(\al+k-1)\om.
\end{equation}
Here we may assume $\al+k-1\ne 0$, since we have $\ddd\xi=0$ if $\xi\in(\Om_{Y,0}^{n'-2})_0$. (Note that $\xi$ may be nonzero if $n'=2$.) So \eqref{18.7} follows from \eqref{18.8}.
\sk
It now remains to show the injectivity of the morphism \eqref{18.2}. For $k\gg 0$, this follows from its surjectivity by using \eqref{16.5}. Indeed, the pole order filtration $P$ is exhaustive, and we have the isomorphisms
\begin{equation} \label{18.9}
H''_{\al+k}=H^{n'-1}(F_h)_{\ee(-\al)}\q\h{if}\,\,\,k\gg 0,
\end{equation}
compatible with the morphism $t:H''_{\al+k}\to H''_{\al+k+1}$ (where $\ee(-\al):=\exp(-2\pi i\al)$), see \cite{Br} (and also Remark~\ref{R18.1} below). This implies the injectivity for any $k$ by using \eqref{18.3} together with the injectivity of the action of $t$ on $H''_h$, see the remark after \eqref{17.3}. This finishes the proof of Proposition~\ref{P18.1}.
\end{proof}

\begin{rem} \label{R18.1}
The action of $\dd_t^{-1}$ on $H''_h$ can be defined in a compatible way with \eqref{17.2}, and the Gauss-Manin system $G_h$ is the localization of the Brieskorn module $H''_h$ by the action of $\dd_t^{-1}$. It is a graded $\C\{t\}\langle\dd_t,\dd_t^{-1}\rangle$-module, and \eqref{17.3} holds with $H''_h$ replaced by $G_h$. We have moreover the canonical isomorphisms
\begin{equation} \label{18.10}
G_{h,\al}=H^{n'-1}(F_h)_{\ee(-\al)},
\end{equation}
compatible with the morphism $\dd_t:G_{h,\al}\to G_{h,\al-1}$.
Note that the condition $N=0$ is equivalent to that $\dd_tt=\al$ on $H''_{h,\al}$. So the actions of $t$ and $\dd_t^{-1}$ on $H''_{h,\al}$ (and also on $G_{h,\al}$ ($\al\ne 0$)) can be identified with each other up to a non-zero constant multiple.
\sk
We have the canonical isomorphisms (see \cite{ScSt}, \cite{Va1}):
\begin{equation} \label{18.11}
F^{n'-1-k}H^{n'-1}(F_h)_{\ee(-\al)}=H''_{h,\al+k}\q\q\bl(\al\in(0,1],\,\,k\in\N\br),
\end{equation}
using the canonical isomorphism \eqref{18.9} together with the action of $t$ or $\dd_t^{-1}$ on $H''_h$ (by using $N=0$ as is explained above).
\end{rem}

\part{Proof of the main theorems} \label{Pa4}
\nin
In this part we prove Theorems~\ref{T2}--\ref{T4} after recalling some basics of Bernstein-Sato polynomials.

\section{Bernstein-Sato polynomials} \label{S19}
Let $f:X\to\C$ be a holomorphic function on a complex manifold. Let $i_f:X\into X\times\C$ be the graph embedding by $f$, and $t$ be the coordinate of $\C$. We have the canonical inclusion
\begin{equation} \label{19.1}
M:=\D_X[s]f^s\into\B_f:=(i_f)_*^{\D}\OO_X,
\end{equation}
where the last term is the direct image of $\OO_X$ as a $\D$-module, and is freely generated by $\OO_X[\dd_t]$ over the canonical generator $\de(f-t)$. Indeed, $f^s$ and $s$ can be identified respectively with $\de(t{-}f)$ and $-\dd_tt$, see \cite{Ma1}, \cite{Ma2}. We have the filtration $V$ of Kashiwara \cite{Ka2} and Malgrange \cite{Ma2} on $\B_f$ along $t=0$ indexed decreasingly by $\Q$ so that $\dd_tt-\al$ is nilpotent on $\Gr_V^{\al}\B_f$. This induces the filtration $V$ on $M$ and on $M/tM$ so that
\begin{equation} \label{19.2}
b_f(-\al)\ne 0\iff\Gr_V^{\al}(M/tM)\ne 0,
\end{equation}
since $b_f(s)$ is the minimal polynomial of the action of $-\dd_tt$ on $M/tM$ by definition. (Here we assume that $b_f(s)$ exists by shrinking $X$ if necessary.) Note that the support of $\Gr_V^{\al}(M/tM)$ is a {\it closed} analytic subset, and the above construction is compatible with the pull-back by smooth morphisms. These together with Remark~\ref{R20.1} below imply the last inclusion in \eqref{1} in the \nameref{intr}.
\sk
Set
$$G_{i\,}\B_f:=t^{-i}M\subset\B_f[t^{-1}]\q\q\q(i\in\Z).$$
We have the canonical inclusion $\B_f\into\B_f[t^{-1}]$, inducing the isomorphisms
\begin{equation} \label{19.3}
\Gr_V^{\al}\B_f\simto\Gr_V^{\al}\bl(\B_f[t^{-1}]\br)\q\q\q(\al>0).
\end{equation}
This gives an increasing filtration $G$ on $\Gr_V^{\al}\B_f$ for $\al>0$ so that
\begin{equation} \label{19.4}
b_f(-\al-i)\ne 0\iff\Gr^G_i\Gr_V^{\al}\B_f\ne 0\q\q\h{for each}\,\,\,\al\in(0,1],\,i\in\N.
\end{equation}
Indeed, by the isomorphisms
\begin{equation} \label{19.5}
t^i:\Gr^G_i\Gr_V^{\al}\bl(\B_f[t^{-1}]\br)\simto\Gr^G_0\Gr_V^{\al+i}\bl(\B_f[t^{-1}]\br)\q\q(\al\in(0,1],\,i\in\Z),
\end{equation}
we get the canonical isomorphisms
\begin{equation} \label{19.6}
\Gr^G_i\Gr_V^{\al}\B_f=\Gr_V^{\al+i}(M/tM)\q\q(\al\in(0,1],\,i\in\Z).
\end{equation}
\sk
For $\al\in(0,1]$, the filtration $G$ on $\Gr_V^{\al}\B_f$ is a finite filtration such that
\begin{equation} \label{19.7}
\Gr^G_i\Gr_V^{\al}\B_f=0\q\h{unless}\,\,\,\,i\in\N.
\end{equation}
This is closely related to the negativity of the roots of $b_f(s)$, see \cite{Ka1}.
\sk
For any $\beta\in\Q$, we have the canonical decomposition
\begin{equation} \label{19.8}
\B_f/V^{\beta}\B_f=\mopl_{\al<\beta}\,\Gr_V^{\al}\B_f.
\end{equation}
by using the action of $s=-\dd_tt$. This is obtained by taking the inductive limit of the decomposition of $V^{\ga}\B_f/V^{\beta}\B_f$ for $\ga\to-\infty$, and is compatible with the canonical surjections
\begin{equation} \label{19.9}
\B_f/V^{\beta}\B_f\onto\B_f/V^{\beta'}\B_f\q\q(\beta>\beta').
\end{equation}
(This means that \eqref{19.9} corresponds via \eqref{19.8} to the canonical surjection associated to shrinking the index set of direct sums.)
\sk
We then get the {\it asymptotic expansion} for any $\xi\in\B_f:$
\begin{equation} \label{19.10}
\xi\sim\msum_{\al\ges\al_{\xi}}\,\xi^{(\al)}\q\h{with}\q \xi^{(\al)}\in\Gr_V^{\al}\B_f,\,\,\al_{\xi}\in\Q.
\end{equation}
This means that the following equality holds for any $\beta>\al_{\xi}$ via the isomorphism \eqref{19.8}:
\begin{equation} \label{19.11}
\xi\,\,\,\h{mod}\,\,\,V^{\beta}\B_f=\mopl_{\al_{\xi}\les\al<\beta}\,\xi^{(\al)}.
\end{equation}
We have to {\it control} this expansion at least in the case $\xi=f^s$ in order to determine $b_f(s)$. Indeed, we have the decomposition of $f^s$ modulo $V^{\beta}\B_f$ in a compatible way with the action of $\D_X[s]$ by using the action of $\C[s]$ for any $\beta>\al_f$ as in \eqref{19.11}, and {\it the direct summands of $\,f^s$ generate $\,\D_X[s]f^s$ modulo $V^{\beta}\B_f$ over $\D_X[s]$}. So the determination of the {\it leading term} of the expansion is {\it not} enough to calculate $b_f(s)$. A similar phenomenon is known for the Gauss-Manin systems in the isolated singularity case, and this is the reason for which the determination of $b_f(s)$ is so complicated. In the weighted homogeneous isolated singularity case, this problem does {\it not} occur for the Gauss-Manin systems. However, it {\it does} occur for $\B_f$ in the non-isolated homogeneous singularity case, since the situation is quite different in the case of $\B_f$.

\section{Remarks related to Bernstein-Sato polynomials} \label{S20}

\begin{rem} \label{R20.1}
In the notation of Section~\ref{S19}, it is well known that $b_f(s)$ depends only on the divisor $D:=f^{-1}(0)$. Indeed, we have a line bundle $L$ over $X$ corresponding to $D$, and the Euler field corresponding to the natural $\C^*$-action on $L$ and the ideal sheaf of the zero-section of $L$ are sufficient to determine the Bernstein-Sato polynomial by using the filtration $V$ of Kashiwara and Malgrange along the zero-section on the direct image of $\OO_X$ as a $\D$-module by the canonical section $i_D:X\into L$ defined by $D$ (where $i_D$ is locally identified with the graph embedding if a local defining function of $D$ is chosen).
\end{rem}

\begin{rem} \label{R20.2}
For any reduced homogeneous polynomial $f$, we have
\begin{equation} \label{20.1}
\R_f^0\subset\ood\,\Z.
\end{equation}
where the left-hand side is the roots of $b_f(s)$ supported at the origin up to sign. Indeed, setting $X^*\defs\C^n\stm\{0\}$ with $j:X^*\into X\defs\C^n$ the canonical inclusion, we have the canonical isomorphism
\begin{equation} \label{20.2}
j^{\D}_!\M\simto j^{\D}_*\M
\end{equation}
for any holonomic subquotient $\D$-module $\M$ of $\Gr_V^{\al}\B_f|_{X^*}$ if $\al\notin\tfrac{1}{d}\one\Z$. Here $j^{\D}_!$, $j^{\D}_*$ denote the functors between the bounded derived category of $\D$-modules with holonomic cohomology $\D$-modules, which can be defined by using \v Cech-type complexes, see for instance \cite{reg}.
\sk
Taking the blow-up $\pi:\Xt\to X$ at the origin, this isomorphism can be reduced to a similar isomorphism for the inclusion $\Xt\stm E\into\Xt$ with $E$ the exceptional divisor (since $\pi$ is proper), and follows from the Verdier-type extension theorem (which involves only the {\it unipotent} monodromy part, see for instance \cite[Proposition 2.8]{mhm}) choosing a local defining function of $E\sst\Xt$. Indeed, $X^*$ is a $\C^*$-bundle over $E\cong\PP^{n-1}$, and $\Gr_V^{\al}\B_f|_{X^*}$ is a {\it monodromical\one} $\D$-module, whose monodromy around $E$ coincides with $T^{-d}$, where $T$ is the monodromy of the nearby cycle functor for $\pi^*\!f$ (since we have locally $\pi^*\!f\eq z_n^dh$ with $z_1,\dots,z_n$ local coordinates of $\Xt$ and $h\ins\C[z_1,\dots,z_{n-1}]$). (Recall that $\Gr_V^{\al}\B_f$ corresponds to the $\lambda$-eigenspace of the monodromy $T$ on the nearby cycle complex $\psi_f\C_X[n{-}1]$ with $\lambda\eq e^{-2\pi i\al}$ if $\al\notin\Z_{\les 0}$.)
\sk
Using an analogue of the theory of intermediate extensions, the above isomorphism implies the isomorphisms of holonomic $\D_X$-modules
\begin{equation} \label{20.3}
j^{\D}_!\M=\Hc^0j^{\D}_!\M\simto j^{\D}_{!*}\M\simto\Hc^0j^{\D}_*\M=j^{\D}_*\M,
\end{equation}
where the middle term is the intermediate extension. Note that $\Hc^kj^{\D}_!\M\eq\Hc^{-k}j^{\D}_*\M\eq0$ for $k\sgt0$ by the construction using \v Cech-type complexes.
\sk
By the adjunction relations, we have for any $\D_{\{0\}}$-module $\M'$
\begin{equation} \label{20.4}
{\rm Ext}^1_{\D_X}(j^{\D}_!\M,i^{\D}_*\M')={\rm Ext}^1_{\D_X}(i^{\D}_*\M',j^{\D}_*\M)=0,
\end{equation}
where $i:\{0\}\into X$ is the canonical inclusion.
\sk
These imply that any holonomic subquotient $\D_X$-module of $\Gr_V^{\al}\B_f$ supported at the origin is a direct factor of $\Gr_V^{\al}\B_f$, and is detectable by applying $i_{\D}^*$ to $\Gr_V^{\al}\B_f$ (which corresponds to taking the Milnor fiber cohomology) assuming $\al\notin\tfrac{1}{d}\one\Z$. We have however $H^{n-1}(\Ff,\C)_{\la}\eq0$ for $\la^d\ne1$ using the geometric monodromy, since $f$ is homogeneous of degree $d$. We thus conclude that there are no holonomic subquotients of $\Gr_V^{\al}\B_f$ supported at the origin if $\al\notin\tfrac{1}{d}\one\Z$. So \eqref{20.1} follows.
\end{rem}

\begin{rem} \label{R20.3}
In the notation of Remark~\ref{R20.2} above, assume the hypotheses~(IS)-(WH) in Theorem~\ref{T2} and moreover $\al\in(0,1)\cap\ood\,\Z$. Then
\begin{equation} \label{20.5}
\aligned\Gr^W_i\Gr^{\al}_V\B_f=0\,\q\h{for}\,\,\,\,|i-n'|>1,\\{\rm supp}\,\Gr^W_i\Gr^{\al}_V\B_f\subset\{0\}\q\h{for}\,\,\,\,|i-n'|=1,\endaligned
\end{equation}
where $W$ is the weight filtration which has symmetry with center $n'\defs n{-}1$.
\sk
If the asymptotic expansion of some $\xi\in F_0\B_f=\OO_X\de(t{-}f)$ in \eqref{19.10} has a nonzero direct factor satisfying
\begin{equation} \label{20.6}
\xi^{(\al)}\in W_{n'-1}\Gr^{\al}_V\B_f\setminus\D_X[s](F_0\Gr^{\al}_V\B_f)\q\bl(\al\in(0,1)\cap\ood\,\Z\br),
\end{equation}
then this would give an example where the implication $\Longrightarrow$ in \eqref{3} fails with condition~\eqref{2} unsatisfied in Theorem~\ref{T1}. Note that the direct image of $\Gr^{\al}_V\B_f$ by the projection to $S=\C$ defined by the function $y$ in Section~\ref{S3} contains direct factors isomorphic to
$$E:=\D_S/\D_S(y\dd_yy),$$
where $y$ is identified with the natural coordinate of $S=\C$. We have the isomorphism
$${\rm DR}(E)_0=C\bl(\Gr_{V_y}\dd_y:\psi_{y,1}E:=\Gr_{V_y}^1E\to\varphi_{y,1}E:=\Gr_{V_y}^0E\br)_0,$$
with $V_y$ the filtration of Kashiwara and Malgrange along $y=0$, and the image of $\Gr_{V_y}\dd_y$ coincides with the subspace obtained by taking $\varphi_{y,1}$ of the intersection of $E$ with the direct image of $W_{n'-1}\Gr^{\al}_V\B_f$.
\end{rem}

\begin{rem} \label{R20.4}
By a calculation of multiplier ideals in \cite[Proposition 1]{bcm} and using \cite{BuSa1}, we get
\begin{equation} \label{20.7}
\dim\Gr_V^{\al}\OO_X=\tbinom{k-1}{n-1}\q\h{for}\,\,\,\al=\kod\in\bl[\nod,\al_Z\br),
\end{equation}
where $\OO_X=\Gr^F_0\B_f$. By \eqref{2.5}, \eqref{2.8}, \eqref{3.8}, this implies
\begin{equation} \label{20.8}
M_k^{(\infty)}\ne 0\q\h{for}\,\,\,\kod\in\bl[\nod,\al_Z\br).
\end{equation}
Indeed, $\Gr_V^{\al}\Gr^F_0\B_f=\Gr^F_0\Gr_V^{\al}\B_f$ is supported at the origin for $\al\in\bl[\nod,\al_Z\br)$, and we have the inclusion $F^p\subset P^p$ ($\forall\,p\in\Z$) together with the vanishing $P^nH^{n-1}(\Ff,\C)=0$. (Note that $\al_Z=\min\bl(\alt_Z,1\br)\les 1$.)
\end{rem}

\begin{rem} \label{R20.5}
Under the assumption {\rm (IS)}, we can show
\begin{equation} \label{20.9}
\al_f:=\min\R_f=\min\bl\{\al_Z,\nod\br\}.
\end{equation}
\sk
In the notation of Section~\ref{S19}, we index the Hodge filtration $F$ on $\B_f$ so that
\begin{equation} \label{20.10}
F_0\B_f=\OO_X\,\de(t{-}f).
\end{equation}
By \cite{BuSa1}, we can essentially identify the induced filtration $V$ on it with the multiplier ideals of $D:=f^{-1}(0)\subset X$ (except the difference between left and right continuities at jumping coefficients). This implies that $\al_f$ coincides with the log canonical threshold, that is, the minimal jumping number, of $D\subset X$. Indeed, we have
\begin{equation} \label{20.11}
\D_X[s]f^s\subset V^{\al}\B_f\q\h{if}\q\OO_X\,\de(t{-}f)\subset V^{\al}\B_f.
\end{equation}
So we can replace $\al_f$, $\al_Z$ in \eqref{20.9} with the log canonical threshold of $D\subset X$ and that of $D\setminus\{0\}\subset X\setminus\{0\}$. The assertion then follows from \cite[Theorem 2.2]{bcm}.
\end{rem}

\section{First step to the proof of Theorem~\ref{T2}} \label{S21}
In this section we prove the following.

\begin{thm} \label{T21.1}
With the notation and assumptions of Section~{\rm\ref{S4}} $($which include {\rm (WH)} in Theorem~$\ref{T2})$, we have the equality
\begin{equation} \label{21.1}
\msum_{k\in\N}\,\dim_{\C}{\rm Coker}\bl(r_k:M_k\to\mopl_{z\in\Si}\,\Xi_{h_z}^{k/d}\br)=\dim_{\C}H^{n'-1}(\Ff,\C),
\end{equation}
where $r_k$ is the composition of
$$M_k\to\Mt_k\q\h{and}\q\Mt_k\to\mopl_{z\in\Si}\,\Xi_{h_z}^{k/d},$$
induced by {\rm\eqref{3.6}} and {\rm \eqref{4.4}--\eqref{4.5}} respectively.
\end{thm}

\begin{proof}
The assertion follows from Theorem~\ref{T10.1} and Proposition~\ref{P18.1} together with \cite{Di1}, \cite{ggr}. In the notation of \cite[Section 1.8]{ggr}, set
$$\Lc_{\la}=\Sc^i\q\q\bl(i=0,\dots,d-1,\,\,\la=\ee(i/d)\br),$$
with $\ee(\al):=\exp(2\pi i\al)$. This is a locally free $\OO_Y$-module of rank 1 having a canonical meromorphic connection so that the localization $\Lc_{\la}(*Z)$ of $\Lc_{\la}$ along $Z$ is a regular holonomic $\D_Y$-module with
$${\rm DR}(\Lc_{\la})=\Rb j_*L_{\la}.$$
The restriction of $\Lc_{\la}$ to the smooth points of $Z$ is the Deligne extension with residue $i/d$, see also \cite[1.4.1]{BuSa2}. On $Y':=Y\setminus\{u=0\}$, there is an isomorphism of regular holonomic $\D_{Y'}$-modules
$$\Lc_{\la}(*Z)|_{Y'}\cong\OO_{Y'}(*Z')h^{i/d}\q\h{inducing}\q\Lc_{\la}|_{Y'}\cong\OO_{Y'}h^{i/d},$$
with $h:=f/u^d$, $Z':=Z\cap Y'$. We have the pole order filtration $P$ on $\Lc(*Z)$ defined by
$$P_j\bl(\Lc_{\la}(*)\br)=\begin{cases}\Lc_{\la}((j+1)Z)&\h{if}\,\,\,j\ges 0,\\ \,0&\h{if}\,\,\,j<0,\end{cases}$$
which induces the pole order filtration $P$ on the Milnor cohomology $H^{n-1}(\Ff)_{\la}$ by taking the induced filtration on the de Rham complex (where the shift of filtration by $n{-}1$ occurs), see \cite[Section 1.8]{ggr}, \cite{Di1}. By using the Bott vanishing theorem, the filtered de Rham complex can be described by using the differential $\ddd_f$ on $\Om^{\ssb}$ defined by
$$\ddd_f\,\om:=f\ddd\om-\kod\,\dfw\om\q\q(\om\in\Om^{\ssb}_k),$$
(see \cite[4.1.1]{Do} and also \cite[Ch.~6, 1.18]{Di1}). We have moreover the surjection
\begin{equation} \label{21.2}
M_{(q+1)d-i}^{(2)}\onto\Gr_P^{n-1-q}H^{n-1}(\Ff,\C)_{\la}\,,
\end{equation}
using the pole order spectral sequence together with \cite[Ch.~6, Theorem 2.9]{Di1}. Here $[\om]\in M_{(q+1)d-i}^{(2)}$ with $\om\in\Om^n_{(q+1)d-i}$ is sent to
$$\bl[i_E(\om/f^{q+1})f^{i/d}\br]\in\Gr_P^{n-1-q}H^{n-1}(U,L_{\la})\q\q\bl(\la=\ee(i/d)\br),$$
with $i_E$ the contraction by the Euler field $E:=\msum_j\,x_i\dd_{x_i}$, using the isomorphism \eqref{9.2}. Here $f^{i/d}$ is a symbol denoting a generator of $\Lc_{\la}(i)\cong\OO_Y$ (which is unique up to a non-zero constant multiple, and is closely related to the above differential $\ddd_f$).
\sk
The target of \eqref{10.2} is calculated by Proposition~\ref{P18.1}. (There is a shift of filtration by 1 coming from \eqref{16.5}.) The graded pieces of the pole order filtration $P$ of restriction morphism \eqref{10.2} is then induced up to a non-zero constant multiple by the morphism $r_k$ in \eqref{21.1}. Here the source of $\Gr_P^{n-1-q}$ of \eqref{10.2} can be replaced by $M_{(q+1)d-i}^{(2)}$ as above, and $r_k$ can be defined by using the contraction with the Euler field together with the substitution $y=1$ after restricting to
$$Y':=\{y\ne 0\}\subset Y=\PP^{n-1}.$$
Note that this is compatible with the blow-up construction in the proof of Theorem~\ref{T5.1}, where the Euler field $E$ becomes the vector field $y\dd_y$ using the product structure of the affine open piece of the blow-up. 
\sk
We have the strict compatibility of \eqref{10.2} with the pole order filtration $P$ by Remark~\ref{R21.1} below. (Here $F$ is the Hodge filtration. We need Proposition~\ref{P18.1} together with \eqref{16.5} and \eqref{18.11} to show $F=P$ on the target.) So $\Gr_P^{\ssb}$ of the cokernel of \eqref{10.2} coincides with the cokernel of $\Gr_P^{\ssb}$ of \eqref{10.2}, and the latter is given by the cokernel of $r_k$ by the above argument. We thus get \eqref{21.1}. This finishes the proof of Theorem~\ref{T21.1}.
\end{proof}

\begin{rem} \label{R21.1}
Let $\phi:(V;F,P)\to(V';F,P)$ be a morphism of bifiltered vector spaces with $F\subset P$ and $F=P$ on $V'$. If $\phi$ is strictly compatible with $F$, then it is so with $P$.
\sk
Indeed, the strict compatibility for $F$ means $\phi(F^pV)=F^pV'\cap\phi(V)$, and
$$\phi(F^pV)\subset \phi(P^pV)\subset P^pV'\cap\phi(V)=F^pV'\cap\phi(V).$$
We thus get $\phi(P^pV)=P^pV'\cap\phi(V)$, that is, $\phi$ is also strictly compatible with $P$.
\end{rem}

\section{Second step to the proof of Theorem~\ref{T2}} \label{S22}
From Theorems~\ref{T5.1} and \ref{T21.1}, we can deduce the following.

\begin{thm} \label{T22.1}
We have
\begin{equation} \label{22.1}
\msum_{q\in\N}\,\dim N_{q+d}\cap\bl(\mopl_{z\in\Si}\,\Xi_{h_z}^{q/d}\br)=\dim H^{n-2}(\Ff,\C),
\end{equation}
where the intersection on the left-hand side is taken in $\Nt_{k+d}$ via the isomorphism {\rm\eqref{4.4}} using the function $y$.
\end{thm}

\begin{proof}
By \eqref{3.6} and \eqref{4.4}--\eqref{4.5}, we have the inclusions
\begin{equation} \label{22.2}
V_1:=N_{q+d}\,\,\,\,\subset\,\,\,\,V:=\Nt_{q+d}\,\,\,\,\supset\,\,\,\,V_2:=\mopl_{z\in\Si}\,\Xi_{h_z}^{q/d},
\end{equation}
where $p+q+d=nd$. By the graded local duality as is explained in \cite[1.1]{kosz} together with the compatibility of duality isomorphisms in Theorem~\ref{T5.1}, we have
\begin{equation} \label{22.3}
V^*=\Mt_p,\q V_1^{\perp}=M''_p,\q V_2^*=\mopl_{z\in\Si}\,\Xi_{h_z}^{p/d},
\end{equation}
where $V^*$ denotes the dual vector space of $V$. Using the diagram of the nine lemma, we get a surjection of short exact sequences
\begin{equation} \label{22.4}
\begin{gathered}
\xymatrix{0\ar[r]&V_1^{\perp}\ar[r]\ar@{->>}[d]&V^*\ar[r]\ar@{->>}[d]&V_1^*\ar[r]\ar@{->>}[d]&0\\ 0\ar[r]&V_1^{\perp}/V_1^{\perp}\cap V_2^{\perp}\ar[r]&V_2^*\ar[r]&(V_1\cap V_2)^*\ar[r]&0}
\end{gathered}
\end{equation}
which implies the exact sequence
\begin{equation} \label{22.5}
V_1^{\perp}\to V_2^*\to(V_1\cap V_2)^*\to 0.
\end{equation}
By Theorem~\ref{T5.1} this gives the exact sequence
\begin{equation} \label{22.6}
M''_p\buildrel{r_p}\over\to\mopl_{z\in\Si}\,\Xi_{h_z}^{p/d}\to\bl(N_{q+d}\cap\bl(\mopl_{z\in\Si}\,\Xi_{h_z}^{q/d}\br)\br)^*\to 0.
\end{equation}
The assertion \eqref{22.1} then follows from Theorem~\ref{T21.1}. This finishes the proof of Theorem~\ref{T22.1}.
\end{proof}

\section{Proof of Theorem~\ref{T2}} \label{S23}
By the pole order spectral sequence, we have the inequality
\begin{equation} \label{23.1}
\msum_{k\in\N}\,\dim{\rm Ker}\bl(\ddd^{(1)}:N_{k+d}\to M_k\br)\ges\dim H^{n-2}(\Ff,\C),
\end{equation}
where the equality holds if and only if the spectral sequence degenerates at $E_2$.
\sk
On the other hand, the proof of \cite[Theorem 5.3]{kosz} actually proves the inequality
\begin{equation} \label{23.2}
\msum_{k\in\N}\,\dim{\rm Ker}\bl(\ddd^{(1)}:N_{k+d}\to M_k\br)\les
\msum_{k\in\N}\,\dim N_{k+d}\cap\bl(\mopl_{z\in\Si}\,\Xi_{h_z}^{k/d}\br).
\end{equation}
So the $E_2$-degeneration follows from \eqref{23.1}--\eqref{23.2} and Theorem~\ref{T22.1}. This finishes the proof of Theorem~\ref{T2}.

\section{Proof of Theorem~\ref{T3}} \label{S24}
The above proof of \eqref{23.2} shows the inequality
\begin{equation} \label{24.1}
\msum_{k\in\N}\,\dim{\rm Ker}\bl(\ddd^{(1)}:N_{k+d}\to M''_k\br)\les
\msum_{k\in\N}\,\dim N_{k+d}\cap\bl(\mopl_{z\in\Si}\,\Xi_{h_z}^{k/d}\br).
\end{equation}
Note that the intersection of $N_{k+d}$ with the kernel of $\ddd^{(1)}:\Nt_{k+d}\to\Mt_k$ is studied in \cite[5.1.2]{kosz}. (The latter morphism is also constructed there.) This implies the injectivity of the composition in Theorem~\ref{T3}, since the failure of this injectivity implies a strict inequality in \eqref{23.2} by using \eqref{24.1}. This finishes the proof of Theorem~\ref{T3}.

\section{Proof of Theorem~\ref{T4}} \label{S25}
We may assume that $d\gess5$ so that $f$ is not ``special" in the sense of Section~\ref{S8} by \eqref{8.6}, since the argument in the case $d\less4$ is trivial (using $\mu_3\eq1$). In the notation of Section~\ref{S8}, the subvariety $\{g\eq0\}\sst\PP^1$ is {\it reduced,} since $\{y\eq0\}\sst\PP^2$ is general. Hence
\begin{equation} \label{25.1}
\dim\bl(R'/(\dd g)\br)_k=\begin{cases}k{+}1&\h{if}\,\,\,0\less k\less d{-}2,\\2d{-}3{-}k&\h{if}\,\,\,d{-}2\less k\less 2d{-}4,\end{cases}
\end{equation}
and 0 otherwise.
By the non-degeneracy of the Grothendieck residue pairing (see \cite[p.\,659]{GH}, there are monomials $m_k\in R'_k$ with $m_kh\notin(\dd g)$ for any $k\in[0,d{-}3]$. (This can be reduced to the case $k\eq d{-}3$.) This implies a certain nontriviality of the image of $m_h$ in \eqref{8.1}.  Using the exact sequence \eqref{8.3}, we then get that
\begin{equation} \label{25.2}
\nu_{k+d+2}-\nu_{k+d+1}\less\dim G^0_{k+d-1}\less k\q\h{if}\,\,\,0\less k\less d{-}3.
\end{equation}
(This holds even in the $\nu_{d+n}$-nonvanishing case if $h\notin(\dd g)$.) On the other hand we have
\begin{equation} \label{25.3}
\mu_{k+2}-\mu_{k+1}\eq k\q\h{if}\,\,\,0\less k\less d{-}1,
\end{equation}
since $\deg\dd_{x_i}f\eq d{-}1$. These imply that
\begin{equation} \label{25.4}
\mu_{k+2}\mi\nu_{k+d+2}>0\q\h{if}\,\,\,1\less k\less d{-}3.
\end{equation}
So the assertion follows using Theorem~\ref{T1}. This finishes the proof of Theorem~\ref{T4}.

\part{Explicit computations}
 \label{Pa5}
\nin
In this part we calculate some examples explicitly.

\section{Example~I} \label{S26}
Set $\mu'_k=\dim M'_k$, $\mu^{(2)}_k=\dim M^{(2)}_k$, etc., and
$$f_1=x^5+y^4z\q\h{with}\q h_1=x^5+y^4,$$
where $n\eq3$, $d\eq5$, $\tau_Z\eq12$, and $\chi(U)\eq1$, see \eqref{15.2}. In this case $\Si$ consists of one point $p:=[0:0:1]\in\PP^2$, and $(Z,p)$ is defined by $h_1$. We have
\begin{equation} \label{26.1}
\begin{array}{rcccccccccccc}
k:&3 &4 &5 &6 &7 &8 &9 &10 &11 &12 &13 &\cdots\\
\ga_k:&1 &3 &6 &10 &12 &12 &10 &6 &3 &1 &\\
\mu_k:&1 &3 &6 &10 &12 &13 &13 &12 &12 &12 &12 &\cdots\\
\nu_k:& & & & & &1 &3 &6 &9 &11 &12 &\cdots\\
\mu_k^{\scriptscriptstyle(2)}:& & & &1 &1 &1 &1 & & & & &\\
\mu''_k:&1 &3 &6 &9 &11 &12 &12 &12 &12 &12 &12 &\cdots\\
\mu'_k:& & & &1 &1 &1 &1 & & & & &\end{array}
\end{equation}
and
\begin{equation} \label{26.2}
\aligned b_{f_1}(s)&=b_{h_1}(s)\,\mprod_{i=6}^9\bl(s+\tfrac{i}{5}\br),\\ \h{with}\q\q b_{h_1}(s)&=(s+1)\,\mprod_{i=1}^4\mprod_{j=1}^3\bl(s+\tfrac{i}{5}+\tfrac{j}{4}\br).\endaligned
\end{equation}
These are done by using respectively Macaulay2 and RISA/ASIR as in \cite{ex}. Note that $\mu^{(2)}_k=\mu_k-\nu_{k+d}$ by \eqref{9} in the \nameref{intr}, and $b_h(s)$ can be calculated also by applying \eqref{7.2}. These two calculations both imply
\begin{equation} \label{26.3}
\R_{f_1}^0=\bl\{\tfrac{6}{5},\tfrac{7}{5},\tfrac{8}{5},\tfrac{9}{5}\br\}.
\end{equation}
This is an example of an extremely degenerated curve (see Remark~\ref{R8.1}), and is a good example for Theorem~\ref{T3} and Corollary~\ref{C3}.

\section{Example~II} \label{S27}
Let
$$f_2=x^5+x^2y^3+y^4z\q\h{with}\q h_2=x^5+x^2y^3+y^4,$$
where $n\eq3$, $d\eq5$, $\tau_Z\eq12$, and $\chi(U)\eq1$, see \eqref{15.2}. Mote that $\Si$ consists of one point $p:=[0:0:1]\in\PP^2$, and $(Z,p)$ is defined by $h_2$. In this case, it is rather surprising that $h_2$ is {\it quasihomogeneous} with $\tau_{h_2}=\mu_{h_2}=12$, see a remark after Theorem~\ref{T2}. (This cannot be generalized to polynomials of higher degrees as far as tried.) We have
\begin{equation} \label{27.1}
\begin{array}{rcccccccccccc}
k:&3 &4 &5 &6 &7 &8 &9 &10 &11 &12 &13 &\cdots\\
\ga_k:&1 &3 &6 &10 &12 &12 &10 &6 &3 &1 &\\
\mu_k:&1 &3 &6 &10 &12 &12 &12 &12 &12 &12 &12 &\cdots\\
\nu_k:& & & & & & &2 &6 &9 &11 &12 &\cdots\\
\mu_k^{\scriptscriptstyle(2)}:&1 &1 & &1 &1 & & & & &\\
\mu'_k:& & & & & & & & &\end{array}
\end{equation}
\begin{equation} \label{27.2}
\aligned b_{f_2}(s)&=b_{h_2}(s)\,\mprod_{i=3}^4\bl(s+\tfrac{i}{5}\br)\cdot\mprod_{i=6}^7\bl(s+\tfrac{i}{5}\br),\\ \h{with}\q\q b_{h_2}(s)&=(s+1)\,\mprod_{i=1}^4\mprod_{j=1}^3\bl(s+\tfrac{i}{5}+\tfrac{j}{4}\br).\endaligned
\end{equation}
\begin{equation} \label{27.3}
\R_{f_2}^0=\bl\{\tfrac{3}{5},\tfrac{4}{5},\tfrac{6}{5},\tfrac{7}{5}\br\}.
\end{equation}
These are compatible with Theorem~\ref{T3} and Corollary~\ref{C3}, since $M'=0$ in this case.

\section{Example~III} \label{S28}
Let
$$f_3=x^5+x^3y^2+y^4z\q\h{with}\q h_3=x^5+x^3y^2+y^4,$$
where $n\eq3$, $d\eq5$, $\tau_Z\eq11$, and $\chi(U)\eq1$, see \eqref{15.2}. As before $\Si$ consists of one point $p:=[0:0:1]\in\PP^2$, and $(Z,p)$ is defined by $h_3$. However, $h_3$ is {\it not} quasihomogeneous with $\mu_{h_3}=12$, $\tau_{h_3}=11$ in this case. We have by RISA/ASIR
\begin{equation} \label{28.1}
\aligned b_{f_3}(s)&=b_{h_3}(s)\,\bl(s+\tfrac{4}{5}\br)\,\mprod_{i=6}^8\bl(s+\tfrac{i}{5}\br),\\ b_{h_3}(s)&=(s+1)\,\mprod_{i=1}^4\mprod_{j=1}^3\bl(s+\tfrac{i}{5}+\tfrac{j}{4}-\de_{i,4}\de_{j,3}\br),\endaligned
\end{equation}
where $\de_{i,k}=1$ if $i=k$, and $0$ otherwise. So we get
\begin{equation} \label{28.3}
\R_{f_3}^0=\bl\{\tfrac{4}{5},\tfrac{6}{5},\tfrac{7}{5},\tfrac{8}{5}\br\}.
\end{equation}
In this case the pole order spectral sequence degenerates at $E_3$, and $\mu^{(3)}_k\ne 0$ if and only if $k\in\{4,6,7,8\}$ by calculations using computer programs based on \cite{DiSt2}, \cite{nwh}. We then get $\R^0_{f_3}$ by applying Theorem~\ref{T1} since $\R_Z\cap\ood\,\Z=\{1\}$ in this case.

\section{Example~IV} \label{S30}
Let
$$f_4=x^4y^2z+z^7\q\h{with}\q h_4=x^4z+z^7,\,\,h'_4=y^2z+z^7,$$
where $n\eq3$, $d\eq7$, $\tau_Z\eq\mu_Z\eq22+8\eq30$, and $\chi(U)\eq1$, see \eqref{15.2}--\eqref{15.6} and A.\ref{A.1} in Appendix below. We have
$$\begin{array}{rccccccccccccccccccc}
k:&3 &4 &5 &6 &7 &8 &9 &10 &11 &12 &13 &14 &15 &16 &17 &18 &19\\
\ga_k:&1 &3 &6 &10 &15 &21 &25 &27 &27 &25 &21 &15 &10 &6 &3 &1\\
\mu_k:&1 &3 &6 &10 &15 &21 &25 &28 &30 &31 &31 &30 &30 &30 &30 &30 &30\\
\nu_k:& & & & & & & &1 &3 &6 &10 &15 &20 &24 &27 &29 &30\\
\mu_k^{\scriptscriptstyle(2)}:& &1 &1 &1 &2 &2 &2 &2 &1 &1 &1 \\
\nu_k^{\scriptscriptstyle(2)}:& & & & & & & & &1 &1 &1 &2 &1 &1 &1 \\
\de_k:& & & & & &1 &1 &1 &1 &1 &1 & & & & & & &\\\end{array}$$
where we have $\de_k:=\mu_k-\nu_{k+d}=\mu'_k$ in this case. The calculation of $\mu_k^{(2)}$, $\nu_k^{(2)}$ is made by computer programs based on \cite{DiSt2}, \cite{nwh}.
\sk
On the other hand we get by A.\ref{A.1} in Appendix below and by RISA/ASIR
$$\aligned 14\,\R_{h'_4}\subset 14\,\R_{h_4}=14\,\R_Z&=\{5\}\sqcup\{7,\dots,21\}\sqcup\{23\},\\
14\,\R_{f_4}&=\{5\}\sqcup\{7,\dots,24\}\sqcup\{26\},\\
\R_{f_4}^0&=\bl\{\tfrac{11}{7},\tfrac{12}{7},\tfrac{13}{7}\br\}.\endaligned$$
Here $3/7\in(\R_Z+\Z_{<0}\bl)\setminus\R_Z$ with $3/7>\al_Z=5/14$, and $\mu_3=\nu_{10}=1$.
(There is a similar phenomenon in the case $f=x^6y^3-z^9$.)

\begin{rem} \label{R30.1}
It is known that for $f\eq x^5\pl y^5\pl yz(x^3\pl z^2w)$, $\Supp\{\mu'_k\}$ is not discretely connected, according to Aldo Conca, see \cite{Sti}. Indeed, $\sum_k\mu'_k\,T^k\eq T^9\pl T^{11}$ in this case (using A.\ref{A.3} in Appendix below). We have $\sum_k\mu'_k\,T^k\eq T^{12}\pl T^{16}$ if $f\eq x^7\pl y^7\pl yz^2(x^4\pl z^3w)$. Such an example is unknown if condition~(WH) is satisfied.
\end{rem}

\part*{Appendix: Explicit calculations using computers} \label{Ap}
\nin
In this Appendix we explain how to use the computer programs Macaulay2 and Singular for explicit computations of the roots of Bernstein-Sato polynomials.

\begin{A}[\bf Computation of roots coming from {\it Z}\one] \label{A.1}
Using \eqref{7.2}--\eqref{7.3} together with the computer program Macaulay2, one can easily calculate the Steenbrink spectrum and the Bernstein-Sato polynomial of a weighted homogeneous polynomial $h$ with weights $w_i$ and having an isolated singularity.
\sk
To determine the weights of weighted homogeneous polynomials with isolated singularities for $n=2$, it is enough to consider the following cases up to a permutation of variables:
$$x^i+y^j,\q x(x^i+y^j),\q xy(x^i+y^j),$$
where the weights $(w_1,w_2)$ are given respectively by
$$\bl(\tfrac{j}{ij},\tfrac{i}{ij}\br),\q\bl(\tfrac{j}{(i+1)j},\tfrac{i}{(i+1)j}\br),\q\bl(\tfrac{j}{ij+i+j},\tfrac{i}{ij+i+j}\br).$$
Combined with \eqref{7.5}, this implies that the Milnor numbers are respectively
$$(i-1)(j-1),\q (i+1)(j-1)+1,\q(i+1)(j+1).$$
\sk
Set $e={\rm GCD}(i,j)$, $a=j/e$, $b=i/e$, and let $m$ be either
$$ij/e\q\h{or}\q(i+1)j/e\q\h{or}\q(ij+i+j)/e,$$
depending on the above three cases. One can calculate the right-hand side of \eqref{7.3} with $T=t^{1/m}$ in the case $h=xy(x^3+y^2)$, for instance, by typing in Macaulay2 as follows:
\ms
\vbox{\small\sf\verb#A=QQ[T]; i=3; j=2; e=gcd(i,j); a=j//e; b=i//e; m=(i*j+i+j)//e#
\sk
\verb#((T^a-T^m)/(1-T^a))*((T^b-T^m)/(1-T^b))#}
\msn
(This can be copied from a pdf file usually.) In this case, the output should be
$$T^{17}\pl T^{15}\pl T^{14}\pl T^{13}\pl T^{12}\pl 2T^{11}\pl T^{10}\pl T^9\pl T^8\pl T^7\pl T^5$$
In the other two cases, the definition of $m$ should be replaced by {\small\sf\verb#i*j//e#} or {\small\sf\verb#(i+1)*j//e#}.
\sk
These calculations can be used to determine $\R_Z$. (For this it would be also possible to use RISA/ASIR if the singularities are not very complicated.)
\end{A}

\begin{A}[\bf Determination of the singularities of {\it Z}\one] \label{A.2}
It is sometimes difficult to see whether a semi-weighted-homogeneous polynomial $h_z$ is quasi-homogeneous or not, see Remark~\ref{R4.1}. In our case this can be verified by seeing whether the total Tjurina number $\tau_Z$ coincides with the total Milnor number $\mu_Z$. The latter can be obtained by using the method in A.\ref{A.1}, and the former is given by high coefficients of the polynomial {\small\sf mu3} in A.\ref{A.3} below. This can be seen explicitly if the separator ``\,;\," at the end of the definition of {\small\sf mu} is removed (where Return must be pressed).
\sk
However, it may be possible that $\tau_Z$ and $\mu_Z$ obtained by the above procedure coincide even though some singularity is {\it not} quasi-homogeneous. This may occur if there is a {\it hidden} singular point of $Z$. It may be verified by typing (or copying) for instance the following in Macaulay2:
\ms
\vbox{\small\sf\verb#R=QQ[x,y,z]; f=x^5*y*z+x^4*y^2*z+z^7;#
\sk
\verb#S=R/(f); I=(radical ideal singularLocus S); decompose I#}
\msn
The output in this case should be
$$\{\h{ideal (z, y), ideal (z, x{+}y), ideal (z, x)}\}.$$
This shows that there is a hidden singular point at $[1:-1:0]$.
\sk
There is another method to see whether all the singular points of $Z$ are quasi-homogeneous by using a computer program like Singular as in \cite{DiSt2}.
Setting $h:=f|_{z=1}$ after a {\it general coordinate change} of $\C^3$ such that ${\rm Sing}\,Z\cap\{z=0\}=\emptyset$, we can compare
$$\tau_h:=\dim\C[x,y]/(h,h_x,h_y)\q\h{and}\q\mu_h:=\dim\C[x,y]/(h^2,h_x,h_y),$$
using Singular (and also \cite{BrSk}), see \cite{DiSt2}. For instance, setting
$$h=x^5(x+1)z+x^4(x+1)^2z+z^7,$$
(which is the restriction of $f$ to $y=x+1$, and defines $Z|_{x\ne y}\subset\C^2$), we can calculate the total Tjurina number of $h$ for the singular points of $Z|_{x\ne y}$ by typing (or copying) the following in the computer program Singular \cite{Sing}:
\msn
\vbox{\small\sf\verb#ring R = 0, (x,z), dp; poly y=x+1;#
\sk
\verb#poly f=x^5*y*z+x^4*y^2*z+z^7;#
\sk
\verb#ideal J=(jacob(f),f); vdim(groebner(J));#}
\msn
We also get the total Milnor number $\mu_h$ of $h$ by replacing {\small\sf\verb#(jacob(f),f)#} with {\small\sf\verb#(jacob(f),f^2)#} (where {\small\sf\verb#f^2#} should be {\small\sf\verb#f^(n-1)#} in general, see \cite{BrSk}).
The problem is then whether we have
$${\rm Sing}\,Z\cap\{x=y\}=\emptyset.$$
This can be verified by comparing $\tau_h$ and $\tau_Z$, where $\tau_Z$ can be obtained by $\mu_k$ for $k\gg 0$, see also the calculation of $\mu'_k$ as in A.\ref{A.3} below (with $f$ replaced by $f$ in this subsection). Note that $\tau_Z$ is also obtained as $\,$tau$\,$ in the Macaulay2 calculation in A.\ref{A.3} below, which can be seen by removing ``;" after the definition of $\,$tau.
\end{A}

\begin{A}[\bf Calculation of the invariants] \label{A.3}
We can get $\sum_k\,\mu'_k\,T^k$ and $\sum_k\,\de_k\,T^k$ for instance in the case of Example~I in Section~\ref{S26} by typing (or copying) the following in Macaulay2:
\ms
\vbox{\small\sf\verb#R=QQ[x,y,z]; n=dim R; f=x^5+y^4*z;#
\sk
\verb#d=first degree f; d1=n*d-n+1; d2=n*d-n; d3=n*d-n+d;#
\sk
\verb#J=ideal(jacobian ideal(f)); A=frac(QQ[T]);#
\sk
\verb#mu1=sub(hilbertSeries(R/J,Order=>d1),A);#
\sk
\verb#mu2=sub(hilbertSeries(R/J,Order=>d2),A);#
\sk
\verb#mu3=sub(hilbertSeries(R/J,Order=>d3),A);#
\sk
\verb#nu=mu1*T^n-((T^d-T)/(T-1))^n; tau=(mu1-mu2)/T^(n*d-n);#
\sk
\verb#mus=tau*(T^(n*d+1)-1)/(T-1)-sub(nu,{T=>1/T})*T^(n*d);#
\sk
\verb#mup=mu1*T^n-mus#
\sk
\verb#nu2=mu3*T^n-((T^d-T)/(T-1))^n; delta=mu2*T^n-nu2/T^d#}
\msn
Here mup, mus mean $\mu'_{\ssb}$, $\mu''_{\ssb}$ respectively. The outputs should be {\it both\one} $T^9{+}T^8{+}T^7{+}T^6$, since this is an example of an extremely degenerated curve, see Remark~\ref{R8.1}. One can see the intermediate results by removing ``\,;\," (and pressing Return). In the four variable case, {\smaller\sf\verb#QQ[x,y,z]#} should be replaced by {\smaller\sf\verb#QQ[x,y,z,w]#}.
\sk
We can apply the above calculation to Walther's example \cite{Wa2} (see also \cite{bha}) where {\smaller\sf\verb#f#} is given by
\ms
\vbox{\footnotesize\sf\verb#x*y*z*(x+3*z)*(x+y+z)*(x+2*y+3*z)*(2*x+y+z)*(2*x+3*y+z)*(2*x+3*y+4*z);#
\sk
\verb#x*y*z*(x+5*z)*(x+y+z)*(x+3*y+5*z)*(2*x+y+z)*(2*x+3*y+z)*(2*x+3*y+4*z);#}
\msn
In this case, $\sum_k\de_k\,T^k$ is a polynomial of degree 16 and 15 respectively. Since $\tfrac{16}{9}\notin\tfrac{1}{3}\,\Z$ and $\R_Z=\{\tfrac{2}{3},1,\tfrac{4}{3}\}$, this implies that $b_f(s)$ is not a combinatorial invariant of a hyperplane arrangement, see also \cite{bha}.  Note that the above first polynomial gives the same $\{\de_k\}$ as Ziegler's (see \cite{Zi}):
\sk
\vbox{\footnotesize\sf\verb#x*y*z*(x+y-z)*(x-y+z)*(2*x-2*y+z)*(2*x-y-2*z)*(2*x+y+z)*(2*x-y-z);#}
\end{A}

{\makeatletter
\def\section#1#2{}%
\def\@startsection#1#2#3#4#5#6{}%
\def\@seccntformat#1{}%
\def\refname{}%
\makeatother

\par\addvspace{3ex plus 1ex minus .2ex}
\noindent{\normalfont\normalsize\bfseries References}
\addcontentsline{toc}{part}{References}  
\par\vskip 2ex plus .2ex

}
\end{document}